\documentclass[12pt]{amsart}
\usepackage[english]{babel}
\title{The ultrafilter: A peerless tool} 

\author{Labib Haddad}
\address{120 rue de Charonne, 75011 Paris, France}
\email{labib.haddad@wanadoo.fr}

\usepackage{amssymb}
\usepackage{amsmath}
\usepackage{endnotes}

\newcommand{\su}{\subsection*}
\newcommand{\head}{\section*}
\newcommand{\noi}{\noindent}

\newcommand{\Ž}{\'e}

\newcommand{\ˆ}{\`a}

\newcommand{\N}{\mathbb N}

\newcommand{\R}{\mathbb R}

\newcommand{\Z}{\mathbb Z}

\newcommand{\cal}{\mathcal}

\newcommand{\leqs}{\leqslant}

\newcommand{\geqs}{\geqslant}

\newcommand{\ali} {\begin{aligned}}   
\newcommand{\ala} {\end{aligned}}

\newcommand{\stm}{\smallsetminus}

\newcommand{\vide}{\emptyset}

\newcommand{\inc}{\subset}

\newcommand{\bc}{\begin{cases}}
\newcommand{\ec}{\end{cases}}

\newcommand{\ba}{\begin{array}}
\newcommand{\ea}{\end{array}}

\let\oldqedhere\qedhere
\renewcommand{\qedhere}{\pushQED{\qed}\oldqedhere}

\begin{document}
\maketitle

\thispagestyle{empty}

\markboth{\sc Labib Haddad}{\sc The ultrafilter: A peerless tool}

\su{Abstract}The ultrafilter, a peerless tool. This paper was meant for a series of talks at the Bratislava Workshop on the density concept, May 2004. A number of the very many facets of ultrafilters are reviewed (some of them, a bit cursorily, as is to be expected in a short space, and time) including Condorcet's Paradox, ultraproducts and the theory of infinitesimals (non-standard analysis), Banach generalized limits in sequence spaces, Choquet's limits for families of closed sets in general topology and intrinsic geometry, representations of topologies as binary relations among ultrafilters, additive bases in number theory. Peerless indeed, and Protean!

\su{Aknowledgement} In his report, the anonymous referee had established a list of questions asking for details concerning \lq\lq certain notions sketched a bit too quickly". He also mentioned about thirty \lq\lq typos" in the first version of the text. I want to thank him most sincerely, for the care he took in his reading. I, of course, corrected those \lq\lq printing errors" and inserted, in braces, the clarifications he asked for. I also want to thank Georges Grekos who read very carefully the first version of the text, drew my attention to typing errors, and provided me with many a wise, and friendly, advice which I have taken in fullest account.

\

Since this text was originally meant to be a sequence of lectures it retains a bit of the {\it oral style} which, I hope, will not affect clarity.

\

The concept of ultrafilter, since that is what it's all about, was introduced by {\bf Henri Cartan} in two short notes in the Comptes Rendus de l'Acad\Žmie des Sciences, Paris, in 1937, for the needs of general topology [{\bf  CARTAN, H.}, Th\Žorie des filtres, C. R. Acad. Sc Paris, 205 (1937) 595-598.; Filtres et ultrafiltres, ibid. 777-779].
Ultrafilters, as is well known, nowadays, \lq\lq are objects situated at the crossroads of almost all chapters of mathematics (set theory, algebra,
topology, analysis, logic, probabilities). Therefore, they are likely to very many and varied presentations, illustrations, and interpretations." [{\bf HADDAD, L.}, Condorcet et les ultrafiltres, {\it in} Math\Žmatiques finitaires et analyse nonstandard, {\it Publ. Math. Univ. Paris VII}, $n^\circ$ 31 (1989) tome 2, p. 343-360. Text of a lecture given at Luminy in 1985].

\

I shall, readily, recall some of their definitions, all equivalent, of course. This should make clear their protean character, whence their wealth.
However, I will start with a small digression that will eventually bring us back, all the same, to our subject.

\

\su{1 An imaginary election} Imagine a triangular election where the electoral college has to choose among three candidates $A, B, C$. Imagine that this electoral college is divided into three trends. One of the three trends prefers candidate $A$ to candidate $B$, and $B$ to $C$. Another has in mind candidate $B$ first, then $C$ followed by $A$. Finally, for the last trend, candidate $C$ is first then $A$ then $B$. Suppose that the supporters of the three trends are equal in number. 

\

The following diagram helps to summarize the situation.

\

$A>B>C \ \ 1/3$ 

$B>C>A \ \ 1/3$ 

$C>A>B \ \ 1/3$.

\

One third of the voters are supporters of $A$, another third are supporters of $B$, and the last third are supporters of $C$. There is no majority to elect one of the candidates at the first round. Imagine any mode of voting you wish, for example the mode for the presidential election in France (which, in this case, takes into account the age of the candidates because of the equality in the distribution of the number of votes). Imagine that one of the three candidates is finally elected.

\

{\it The day after election.} In any such situation, there will always be a majority of disgruntled. The phenomenon is not so rare and has already been noticed before, since the day men exist ... and vote. However, there is a remarkable thing about the present case. Whatever the method of voting is, whoever the winner is, there will always be a majority which prefers another candidate. For example, if $A$ is designated, there is a comfortable two-thirds majority that would prefer $C$ to be elected. Similarly, if $C$ was elected, there will still be a two-thirds majority that prefers $B$. If $B$ had won the election, a two-thirds majority would prefer $A$.

\

{\bf In a triangular election, there might exist a majority of voters who rather than the winner himself prefer one of the defeated candidates.}

\

That is, briefly presented, one of the forms of the Condorcet paradox.

\

{\it The day after election}. Even if the voters give up trying to designate a winner, if they are content with establishing just a simple ranking, this would still lead to an equally paradoxical situation where $A$ is preferred to $B$, $B$ to $C$ and $C $ to $A$, each of these preferences being the wish of a two-thirds majority, never the same, each different! A situation \lq {\it \ˆ la Condorcet}'.

\

{\bf Attempting to aggregate individual preferences might lead to an \lq\lq inconsistent" collective circular ranking.}

\

Different situations and percentages can still produce a similar result: A majority of disgruntled, whoever the winner, ready to have one of the unsuccessful candidates elected, as the following example shows.

\

Three candidates,  $A, B, C$. Among a hundred of voters, 45 are in favor of $A$, 25 in favor of $B$, and 30 in favor of $C$. Among the 45 supporters of $A$, there are 20 who prefer $B$ to $C$, and 25 who prefer $C$ to $B$. The 25 supporters of $B$ divide as follows: 10 prefer $A$ to $C$ and 15 prefer $C$ to $A$. Finally, the 30 supporters of $C$, all, prefer $B$ to $A$. This is summarized in the following sketch:

\

$A>B>C \ \ 20 \ \text{voters}$

$A>C>B \ \ 25$

$B>A>C \ \ 10$

$B>C>A \ \ 15$

 $C>B>A \ \ 30$.
 
 \
 
Let us say it again. Whatever the method of voting is, whoever the chosen candidate is, there would be a majority of disgruntled, but, which is more, a majority would be  ready to overthrow the winner and vote for one of the beaten candidates.

\

{\bf Indeed ...}

\

{\it The day after election.} If $A$ were the winner, 55 voters against 45 would still prefer $B$ to $A$. Similarly, if $B$ is the winner, 55 voters against 45 prefer $C$. Finally, in the case where $C$ is elected, we would end up with 55 against 45 voters preferring $A$.

\

In such a situation, regardless of the procedure, trying to break the deadlock, there will always be a majority ready to agree on another solution.

\

{\it Aftermath of elections}!

\

[{\bf HADDAD, L.}, \'Elections, ultrafiltres, infinit\Žsimaux ou le paradoxe de Condorcet, in Condorcet, Math\Žmaticien, \'Economiste, Philo\-sophe, Homme politique, (p. 87-91), Colloque international, Paris, juin 1988, \Ždition Minerve].

\su{2 Two steps}

\

\

{\bf 1785.}  {\bf Condorcet} was one of the first to highlight and analyze this inevitable consequence of majority voting, in the Discours pr\Žliminaire to his {\it Essai sur l'application de l'analyse \ˆ la probabilit\Ž des d\Žcisions rendues \ˆ la pluralit\Ž des voix}, [Imprimerie Royale, Paris, 1785. (A photographic reprint has been published by Chelsea Publishing Company, 1972. The Discours has also been inserted in the following book: Condorcet, Sur les \Žlections et autres textes, Corpus des \oe uvres de philosophie en langue fran\c caise, Fayard 1986)].

\

{\bf 1952.} {\bf G. Th. Guilbaud} attracting attention again to this phenomenon, named it {\bf Condorcet effect}, in a remarkable and important text [Les th\Žories de l'int\Žr\t g\Žn\Žral et le probl\me logique de l'agr\Žgation, {\it \' Economie appliqu\Že}, {\bf 5} (1952) $n^\circ$ 4, oct.-d\Žc., 501-551. (The paper was printed again as chapter II in the following book: {\it \'El\Žments de la th\Žorie math\Žmatiques des jeux}, Monographies de re\-cherches op\Žrationnelles, 9, Collection directed by G. Morlat, AFIRO, Dunod, Paris, 1968)].

\

{\bf Condorcet effect.} The harm is profound. For example, a club can very well want its emblem to be, all at a time, and by majority votes,
a convex quadrilateral with four equal angles and four equal sides but which would not be a square!

\

If we want to avoid the Condorcet effect, we cannot maintain the majority system. Or, else, we should have to modify it by extending the concept of majority.

\

Guilbaud introduces the notion of {\it a generalized voting system} and analyzes the way it operates.

\

Among all the possible combinations within an assembly $A$, we could, for example, distinguish a number of them, in advance, to hold the role of majorities. We would call them {\it the efficient coalitions} (to distinguish them from majorities in the ordinary sense), and it would be understood that the assembly confirms all the choices (and only those choices) that would be approved by one (any) of those efficient coalitions.

\

\su{3 Voting systems}
Let $A$ be an assembly. Each subset $K$ of $A$ will be called {\bf a coalition}. We let $K^c = A \stm K$ denote the opposing (or complementary) coalition to $K$.

\

A {\bf voting system} for the assembly $A$ is a given set $\cal E$ of coalitions of $A$ (to be used as substitutes for ordinary majorities). In other words, $\cal E$ is a subset of the set $\cal P(A)$ of all subsets of the set $A$. {The elements of $\cal E$ are, by definition, the efficient coalitions.}

\

The system ($A,\cal E$) operates in the following way. Questions submitted to the assembly are dichotomous choices which the members of the assembly answer by {\bf yes} or {\bf no}. Thus two coalitions form, the {\bf pros} and the {\bf cons}. The two coalitions are complementary. [The choice being dichotomous, there is no room for abstention. We will say later why this is not a real constraint since any multiple choice boils down to a suitable set of dichotomous choices.] Of these two coalitions the efficient one outweighs the other, i,e., the one which belongs to $\cal E$.

\

If we, then, want the assembly not to be able, at the same time and with a same impulse, accept and reject a given choice, and if we also want the assembly to be determined in all circumstances, it is necessary (and sufficient) to require the following condition on the set $\cal E$:

\

{\bf C1} A coalition is efficient if and only if the opposing coalition is not.

\

If, moreover, we want that, given two consecutive choices, the assembly will not recant, that is, if we want that, given two incompatible choices, the assembly will not adopt both, it is necessary and sufficient to require, also, the following condition:

\

{\bf C2} Any coalition that contains an efficient coalition is efficient, itself.

\

To say that two choices $u$ and $v$ are incompatible amounts to say that one cannot be {\bf pro} to both $u$ and $v$. In other words, the coalition $u^+$ of the {\bf pro} $u$ is contained in the coalition $v^-$ of the {\bf  con} to $v$, and, similarly, $v^+$ is contained in $u^-$. Abbreviating, both conditions can be written as follows:

\

{\bf C1} \ $K\in \cal E$ if and only if $K^c\notin \cal E$. 

{\bf C2} \ If $L\supset K\in \cal E$ then $L\in \cal E$.

\

We have, thus, made sure that these two conditions are necessary and sufficient for the system to avoid {\it the most elementary failures}: The assembly can always decide, in all circumstances, and each of its decisions cannot contradict any another.

\

{\bf In what follows, all voting systems under consideration are supposed to meet both conditions {\bf C1} and {\bf C2}.}

\

But this in not sufficient to avoid {\bf \ˆ la Condorcet situations}: This is clearly shown by the example of the classical majority system! Indeed, this system meets both conditions {\bf C1} and {\bf C2}, and yet, as we have seen, it is not safe from paradoxical Condorcet situations.

\

\su{4 Remarks}

\

\

It will be observed that condition {\bf C1} has the following consequence: Two distinct members, $x$ and $y$, of the assembly cannot both belong to all the efficient coalitions, since out of the singleton $\{x\}$ and its complement, only one can be efficient!

\

This same condition {\bf C1} also implies that there are as many efficient coalitions as inefficient coalitions.

\

It is also clear that conditions {\bf C1} and {\bf C2} together imply that unanimity is efficient! And the empty coalition is inefficient!

\su{5 Two \lq\lq classical" examples}

\

\

The typical example of voting systems is the majority system, of course, whose efficient coalitions are, by definition, the majority coalitions, i.e., coalitions which outnumber their opposition. This system satisfies condition {\bf C2} but it only satisfies condition {\bf C1} in assemblies having an odd number of members. When an assembly has an even number of members, there is a way very commonly used to avoid the difficulty by introducing  {\it  a casting} vote.

\

The \lq\lq majority system with a {\bf casting vote} for the president" is defined by distinguishing a particular member of the assembly (the \lq\lq Chairman", the \lq\lq oldest", the \lq\lq dean", or any other person designated in advance) and by declaring \lq\lq majority" any coalition with more members than the opposite coalition and, also, any coalition that includes the president and has as many members as its opposite. This system still satisfies both conditions {\bf C1} and {\bf C2}, whether the assembly is even or odd.

\

Notice, by the way, that nothing forbids to consider a \lq\lq majority system with {\bf a minored vote} for a mock-president" as would be the case, for example, with a group of schoolchildren wanting to lower the voice of the youngest!

\

The other best-known system is the {\bf dictatorial} system: A particular member, $d$, of the assembly is designated (or, more often, appoints himself!) Efficient coalitions are those to which the dictator $d$ belongs. 

\

In this dictatorial system, which obviously satisfies conditions {\bf C1} and {\bf C2}, only the opinion of the dictator is taken into account.

\

Notice that singleton $\{d\}$ is then the smallest of all efficient coalitions, i.e., it is contained in each other efficient coalition. Similarly, it will be noticed that, conversely, a voting system satisfying conditions {\bf C1} and {\bf C2} and having a smallest efficient coalition is necessarily dictatorial. Indeed, as has been already noticed (in {\bf 4}), this smallest efficient coalition is necessarily a singleton.

\

\su{6 Some complements}

\

\

\begin{enumerate}

\item   {\bf The case of multiple choice}. - In a multiple choice question, one can encode the set of choices using the binary $\{0, 1\}$-code. This turns the multiple choice into a sequence of dichotomous choices.

\

Doing that, in a vote by \lq\lq yes" or \lq\lq no", we can always make room for, and authorize, abstention as well as blank votes, refusal to vote, or any other modality the community would wish to introduce, without removing anything out from the dichotomous mode of operation described above.

\

\item {\bf The method of weightings}.- There is a variant to the majority system, that of weighting votes. It generalizes, at the same time, the \lq\lq majority system with a casting vote for the president" and the \lq\lq majority system with a minored vote for a mock-president" (as well as the dictatorial system!).

\

To each member $x$ of the assembly $A$  a \lq\lq weight" $p(x)$ is assigned. The weight $p(K)$ of a coalition $K$ is defined to be the sum of the weights of its members. A coalition $K$ is then termed efficient when its weight exceeds that of the opposing coalition $K^c$, i.e., when $p(K) > p(K^c)$. The weights $p(x)$ need not be integers; they can be real numbers [or, even, elements of any given linearly ordered group].

\

It is easily seen that conditions {\bf C1} and {\bf C2} hold in the system thus constructed if (and only if) the weights $p(x)$ are never negative and no coalition has the same weight as its opposite.

\

\item {\bf A counter-example}.- Here is a very simple example of a voting system $(A,\cal E)$ satisfying conditions {\bf C1} and {\bf C2} which cannot be obtained by any weighting. On the finite 7-points projective plane $A$, we define $\cal E$ to be the set of all subsets of the projective plane which contain at least 5 points, and those subsets which contain at least one of the 7 lines in the plane.

\

\item A relatively simple computation shows that, for a three-person assembly, assuming the equidistribution of opinions, the probability of meeting a Condorcet situation, in an election among three candidates, is equal to 1/18, slightly more than 5.5 per-cent of the cases.

\

\item {\bf Three friends}.- Three friends $x,y,z,$ decide to go to the cinema. They have the choice between three films, $a,b,c$. Their preferences are as follows: $x$ would rather see $a$ than $b$ and $b$ rather than $c$; to $y$, $b$ is better than $c$, and $c$ better than $a$; to $z$, $c$ is better than $a$, and $a$ better than $b$. Voting, in order to compare the films pairwise, the three friends rapidly end up with the following circular ranking: $a$ more than $b$ more than $c$ more than $a$. {\bf A beautiful miniature of the Condorcet effect.}

\

Fortunately, the situation here is not so serious. We know that the more persuasive, the most charming of the three friends will win. If this circular ranking does not occur more than 5.5 percent of the cases, no one will worry. If the frequency is still substantially lower, this means that friendship here is coupled with common inclinations. If it is much larger, then we clearly face very strong feelings.

\end{enumerate}

\

\su{7 How to avoid the Condorcet paradox}

\

\

We already said that the majority system can produce a \lq\lq Condorcet effect." Let us observe that only three persons are enough to have a Condorcet effect but that three questions at least are needed.

\

For a voting system $(A,\cal E)$, condition {\bf C1} ensures the consistency of collective decision when a single question is asked about a single binary choice. Condition {\bf  C2}, together with {\bf C1}, ensures collective consistency when two questions are asked. But what must be added to ensure consistency when three (or more) questions are asked, i.e., how to avoid Condorcet's paradox?

\

The answer is simple: It is necessary and sufficient that the intersection of any two efficient coalitions be efficient, itself.

\su{8 Lemma [Guilbaud]}

\

\

{\it A voting system $(A,\cal E)$ (satisfying conditions {\bf C1} and {\bf C2 )} avoids the Condorcet effect, if and only if it has the following property}:

\

{\bf C3} The intersection $K\cap L$ of any two efficient coalitions $K$ and $L$ is also efficient.

\

{\bf  Indeed}, if coalition $H = K\cap L $ were not efficient, we would be in a well-known situation in game theory. One could imagine a sequence of three questions $p,q,r,$ such that a positive answer to the first two implies a positive answer to the third, and also imagine that we have $p^+ = K$, $q^+ = L$, and $r^- = H$. Three questions of the following kind:

\

$p$: \ Do you want a member of coalition $K$ to win the jackpot? 

$q$: \ Do you want a member of coalition $L$ to win the jackpot? 

$r$: \ Do you want a member of coalition $H$ to win the jackpot?

\

\noi So, the assembly would have decided, through coalition $K$, that the jackpot must go to a member of $K$; it would have decided, similarly, through coalition $L$, that the jackpot must go to a member of $L$; it would have finally decided, through coalition $H^c$, that the jackpot should not go to a member of $H$!

\

\noi Conversely, if condition {\bf C3} is satisfied, no inconsistent decision is possible since there will always exist at least one member of the assembly who belongs to all the coalitions that have enforced these decisions and therefore he would have endorsed them all (this member is supposed, by hypothesis, to be himself consistent). \qed

\

Alas! The above proof lets us foresee that the remedy leads to ill...

\

\su{9 Guilbaud's Theorem}

\

\

{\it A voting system meeting all three conditions, {\bf C1}, {\bf C2}, and {\bf C3}, is dictatorial.}

\

{\bf Indeed}, let $M$ be the intersection of all efficient coalitions. Using condition {\bf C3}, and induction, coalition $M$ is seen to be efficient. It, therefore, is the smallest of all efficient coalitions! \qed

\

\su{10 The dead-lock and the way out}

\

\

We are deadlocked, having to choose either Condorcet effect or dictatorship. To get out of the dead-end, one could work out, and this has been done, all kinds of systems less \lq\lq rudimentary" than those described above. But Guilbaud's Theorem has a character of robustness. It also applies to these more complex decision making systems.

\

\noi In 1951, {\bf Kenneth Arrow} [Social choice and individual values, John Wiley \& Sons, New York, 1963] announced his often quoted result on the impossibility for a non dictatorial system to reconcile certain (natural) conditions for consistency. Arrow establishes a general theorem. Guilbaud's Theorem is a special case, an {\bf exemplary} case.

\

\noi \{The theorem of Arrow says that the sole decision making system which meets a list of \lq\lq natural" broad enough compatibility and consistency conditions [too long to enumerate] is the dictatorial system. One can verify, minutely, that conditions {\bf C1}, {\bf C2}, and {\bf C3}, imply, one by one, each of the conditions in the Arrow list. It follows that Guilbaud's Theorem is a special case of Arrow's.\}

\

\noi Either Condorcet effect or dictatorship. Indeed, in one word, the difficulty is such that, no doubt, it will not be enough to blow on it to make it fall. So, was it not our intent sketching thus the route between Condorcet and Arrow.

\

\noi The knot is tight. It deserves better than amused attention. It involves, all together, politics! But also sociology, psychology, economics, and ... mathematics!

\

What is to be done? against the theorem of Arrow. Ah! We had almost forgot to say that Guilbaud's Theorem holds provided  the assembly A is {\bf finite}!

\

\su{11 Ultrafilters}

\

\

Those who know will already have recognized in the set of three conditions, ({\bf C1, C2, C3,}) one of the many equivalent definitions of ultrafilters.

\

\noi One of the other possible definitions, which is particularly adapted to the circumstances, is as follows.

\

\noi An ultrafilter on a set $A$ is a non-empty set $\cal U$ of non-empty subsets of the set $A$ which has the following two properties.

\

\noi For each subset $K$ and $L$ of $A$, we have

\

{\bf U1} \  \  $K\cap L\in \cal U$ if and only if $K\in\cal U$ and $L\in \cal U$, 

{\bf U2} \ \ $K\cup L\in \cal U$ if and only if $K\in \cal U$ or $L\in \cal U$.

\

\noi Notice in these two statements, the link between the operation \lq\lq intersection" with the conjunction \lq\lq and", on the one hand, and of the operation \lq\lq union" with the disjunction \lq\lq or", on the other, which should not unduly surprise!

\

\noi One easily verifies that the set of two conditions ({\bf U1}, {\bf U2},) is equivalent to the set of three conditions ({\bf C1, C2, C3}). In other words, a system of voting $(A,\cal E)$ satisfies the conditions {\bf C1, C2, C3,} if and only if the set of efficient coalitions $\cal E$ is an ultrafilter on the set $A$.

\

\noi We are back to our concern and this provides one more suit to ultrafilters. Ultrafilters considered as voting systems ensuring the con- sistency of decisions, that is, systems of vote avoiding Condorcet effect.

\

\noi To each element $x$ of the set $A$ corresponds the ultrafilter $\cal U_x$  consisting of the subsets of $A$ to which the element $x$ belongs. These are the {\bf trivial} ultrafilters on $A$. They correspond to the dictatorial voting systems.

\

\noi The Theorem of Guilbaud simply reflects the well-known following fact: On finite sets, all ultrafilters are trivial.

\

\noi However, on an infinite set, there always exist non-trivial ultrafilters, provided we use the axiom of choice or another additional (suitable) axiom added to the classical Zermelo-Fraenkel theory ZF.

\

{\bf \lq\lq Truth in mathematics cannot be decided with votes."}
This short sentence is often heard, launched during heated discussions, at the moment the debate begins to produce its most specious arguments.

\

\noi Well! Try all the same. In other words, let us consider a voting system ($I,\cal U$) with $I$ any set and $\cal U$ an ultrafilter on $I$.

\

\su{12 Ultrapowers}

\

\

Let us also consider a set $X$ and ask that every member $i$ of the assembly $I$ chooses an element $x_i$ of the set $X$. We, thus, obtain a family, $x = (x_i)_{i\in I} = (x_i)$, of elements of $X$ indexed by $I$ and which, in some way, is the collection of individual choices of the members of the assembly. But, then, what is the collective choice of the assembly (according to the voting system defined by the ultrafilter $\cal U$)?

\

\noi If there exists a well determined element $e$ of $X$ which an efficient coalition has chosen, there is no hesitation saying that this element represents the collective choice of the assembly since the element was elected by a vote. In other words, if the set $\{i\in I : x_i =e\}$ belongs to the ultrafilter $\cal U$, it is legitimate to say that $x \equiv e \  (\text
{modulo} \ \cal U)$.

\

\noi More generally, if the assembly makes two consecutive choices, let us say, $x = (x_i)$ and $y = (y_i)$, we consider that these two choices are the same when there is an efficient coalition to say so, that is, when the coalition $\{i \in  I : x_i =y_i\}$ belongs to $\cal U$. We shall then say that we have $x\equiv y \ (\text{modulo} \ \cal U)$. This is an equivalence relation on the Cartesian product $X^I$.

\

\noi One might have recognized, straightforwardly, with no detour, the notion of {\bf ultrapower}. Indeed, the set of all \lq\lq collective choices" of the assembly, thus defined, is nothing else but the ultrapower $X^{\cal U}$ of $X$ with respect to the ultrafilter $\cal U$, i.e., the quotient (the set of cosets) of the Cartesian product $X^I$ by the equivalence relation.

\

\noi It is a mathematical construction that can metaphorically be interpreted as the result of the works of an assembly voting according to a given ultrafilter.

\

\noi I have already had the opportunity to develop this {\it metaphor} elsewhere (Luminy 1985, see Condorcet et les ultrafiltres, loc. cit.)

\

\noi One usually associates the name of {\bf \L o\'s} and the year1955, even 1949, to the notion of ultrapower [{\bf \L o\'s, J.}, Quelques remarques,th\Žor\mes et probl\mes sur les classes d\Žfinissables d'alg\bres, in {\it Mathematical interpretation of formal systems}, Amsterdam, 1955; p. 98-113; O matrycach logicznych, {\it Prace Wroclawskiego Towarzystwa Naukowego}, Wroclaw, 1949].

\

\noi We notice by the way that nobody would have missed to notice that this notion is, so to speak, in germ in the notion ... of {\bf germs of functions}, well-known to geometers and analysts.

\

\{Two functions $f :X\to Y$ and $g :X\to Y$ are said to be {\it equivalent} modulo a given filter $\cal F$ on $X$ when there is some $E \in \cal F$ on which $f$ and $g$ coincide. Germs of functions according to the filter $\cal F$ are, by definition, the equivalence classes modulo this filter.\}

\

\noi What happens, however, if, for a given collection $x = (x_i)$ of individual choices, there is no well determined element $e$ of $X$ for which we have $x\equiv e \ (\text{modulo} \  \cal U)$? And this  is inevitably likely to occur when the ultrafilter $\cal U$ is not trivial, and the set $X$ is infinite. Well! We simply say that $x$ (modulo $\cal U$) is the collective choice of the assembly. We would have thus created an ad hoc unreal object, the price to pay in order to avoid, at the same time, dictatorship and inconsistency.

\

\noi {\bf Bourbaki} uses a similar construction in his \lq\lq Structures fondamentales de l'analyse", livre IV, entitled \lq\lq Fonctions d'une variable r\Želle". It serves him to define \lq\lq scales of comparison."

\

\noi In a footnote on page 57 of Chapter 5, in the 1951 edition (the note apparently disappeared from the last \lq\lq gray"  edition) Bourbaki specifies the following about his construction:

\

\lq\lq This calculus on equivalence classes of \lq locally equal' functions plays an important role in many theories that will be developed in this treatise, especially in the theory of differentiable manifolds."

\

\noi An important role! Did he say so? The development of Bourbaki culminates in an Appendix on Hardy fields [{\bf BOURBAKI}, Fonctions d'une variable r\Želle, FVR V.36] which, as will later be noticed, are a first sketch of renovated infinitesimals.

\

\noi A first sketch, because it does not contain the essential remarks about the {\bf permanence of first order formulas} which is {\it the salt of ultraproducts}.

\

\su{13 Ultraproducts}

\

\

More generally still, if each individual $i$ chooses a set $X_i$, to start with, then limits his choices to elements of this set, the collective choices are represented (modulo $\cal U$) by families $(x_i)_{i\in I}$  indexed by $I$ where $x_i$ is an element  of $X_i$ for each individual $i$ in $I$.

\

\noi This is nothing else but the {\bf ultraproduct} of the family of sets $(X_i)$ with respect to $\cal U$. [Of course, each ultrapower is a special case of ultraproducts.]

\

\noi By its votes, the assembly will be able to handle all kinds of mathematical objects and form ultraproducts of groups, of fields, of Banach spaces, of analytic functions, what else do I know ...

\

\noi Among the collective choices of the assembly, some will be those that are conventional, \lq\lq real", and those which are fictional, \lq\lq ways of talking" (to silence dissent!) But with these ways of talking, we can rightly discourse and try to tell truth from falsehood.

\

\noi In order to know whether a proposition is true or false, we will let the assembly vote with respect to the system $\cal U$. And this method succeeds, beyond all hope.

\

\su{14 The language}

\

\

Without going into details, let us say that those who practice  \lq\lq mathematical logic" have developed, worked out, and use \lq\lq formal" languages with many symbols and very strict syntax to be safe from all past, and to come, disputes.

\

\noi The mathematician could well do without all this complexity. He only needs to know that these languages exist and have a knowledge of the basics. Not to shrink in front of this \lq\lq artificial barrier", moving forward, learning how to avoid \lq\lq missteps", not to be impressed or put off by complicated words that cover simple concepts, this is possible and even desirable if one does not want to get lost in the meanders and quicksands of the beautiful mathematical landscape.

\

\noi To be more precise, the proposals that will be submitted to the assembly, are stated, expressed, using a language called \lq\lq first order with equality", and with the binary predicate of membership in a set theory.
In order to talk to the assembly, just use negation, {\bf  \lq\lq not"}, disjunction, {\bf \lq\lq or"}, the phrase {\bf \lq\lq there exists"}, the equality relation, {\bf \lq\lq ="} and the membership relation, $\lq\lq \in"$ [a total of two connectors, one quantifier, and two specific binary predicates]. We will derive all other connectors and quantifiers, including the conjunction {\bf \lq\lq and"}, implication, the quantifier {\bf \lq\lq for all"},  in the most simple, usual way.

\

\noi Thus, for example, \lq\lq $p$ and $q$" means \lq\lq not((not $p$) or (not $q$))". 

\noi Similarly, \lq\lq for all $x$, $p(x)$" means \lq\lq not (exists $x$, (not $p(x)$)".

\

\noi The assembly will understand perfectly well when told \lq\lq for all $x$, $x$ does not belong to $X$" that the $X$ we are talking about is the empty set!

\

\su{15 Truth according to $\cal U$}

\

\

Let $P(x,\dots,y)$ be a statement in the agreed language, submitted to the assembly, where $x = (x_i), \dots , y = (y_i)$ represent the choices of the assembly. Denote by $V(P;x,\dots,y)$, or more simply $V(P)$ when confusion is minimal, the coalition of members who believe that this statement is satisfactory \lq\lq from their own point of view."

\

\noi That is, let $V (P ) = V (P ; x, \dots , y) = \{i \in I : P (x_i, \dots, y_i) \ \text{is true}\}$. This coalition could be efficient or inefficient, i.e., belong or not to the ultrafilter $\cal U$ . If it is efficient, we say that $P (x, \dots , y)$ is true according to $\cal U$. Otherwise, we say that $P(x,\dots,y)$ is false according to $\cal U$.

\

\noi When the statement $P(x,y)$ reads \lq\lq $x = y$", one will thus obtain $V (P ) = \{i \in I : x_i = y_i \}$. This is a familiar coalition already encountered. It follows that the following two expressions are synonymous:

\

(1) $x \equiv y \ (\text{modulo}\ \cal U)$ and

(2) the statement \lq\lq $x = y$" is true according to $\cal U$.

\

\noi They both have exactly the same meaning, and we have to say that everything has  been  done for that purpose!

\

\noi Similarly, for the proposition $Q(x,X)$ which reads \lq\lq $x \in X$", we have $V (Q) = \{i \in I : x_i \in  X_i\}$.
This proposition $Q(x,X)$ is true according to $\cal U$ if and only if a \lq\lq majority", sorry, I mean an efficient coalition of members $i$ choose $x_i \in X_i$.

\

\noi When this is the case, we will still say that $x \in X\ \ (\text{modulo} \  \cal  U)$. 

\

\noi We have thus given  precise meanings to  predicates \lq\lq  = " and \lq\lq$ \in$" in the fictitious world of collective choices, meanings that largely overflow the usual sense. They only overflow since they coincide with the usual sense for all real objects!

\

\noi  What is the relation between truth according to $\cal U$ and truth in our fictitious world? They are one and the same. That is the main result, the most important.

\

\su{16 A Fundamental Lemma} {\it Statement $P (x, \dots , y)$ is true if and only if an efficient coalition so decides.}

\

{\bf In other words, $P(x,\dots,y)$ is true in the fictitious world if and only if $P(x,\dots,y)$ is true according to $\cal U$.}

\

This lemma barely deserves a proof, the language being itself so convincing! The proof is very simple and is by induction \lq\lq on the complexity" of statements. It could make an excellent exercise. It is to be found in all textbooks that consider ultraproducts. [See, for instance, Condorcet et les ultrafilters, loc. cit.].

\

\noi {\bf A sketch of the proof}.-

\

\

1) If statement $P$ only contains one of the two symbols \lq\lq =" or \lq\lq$\in$" (but no logical symbol) the result follows from the definition of the meanings of these two predicates in the fictitious world.

\

2) If statement $P$ is of the form (not $Q$): Suppose the lemma has been established for $Q$ and just notice that $V (\text{not}\ Q)$ and $V (Q)$ are two complementary coalitions!

\

3) If statement $P$ is of the form ($R$ or $S$): Suppose, similarly, that the lemma has already been established for each one of the two statements $R$ and $S$, then notice that $V(R \ \text{or}\ S)$ is the union of $V(R)$ and $V(S)$.

\

4) Finally, if $P(x,\dots,y)$  is of the form \lq\lq exits $z, T(z,x,\dots,y)$", we proceed as follows.

\

\noi For each collective choice $c = (c_i)$ the coalition $W = V(P;x,\dots,y)$  contains $V(T;c,x,\dots,y)$. Indeed, if $i$ belongs to $V(T;c,x,\dots,y)$ then, from $i$'s point of view, the statement $T(c,x,\dots,y)$ is satisfactory, in other words, the statement \lq\lq exists $z, T (z, x_i, \dots , y_i)$ " is true, so that statement $P(x_i,\dots,y_i)$ is also true, so that $i$ also belongs to $W$.

\

\noi Suppose then that the lemma were established for statement $T$ (and whatever the choices involved). Then we reason as follows.

\

\noi 1. Statement $P (x, \dots, y)$ is true if and only if there exists a collective choice $c = (c_i)$ such that $T (c, x, \dots, y)$ is true, which implies that $V(c,x,\dots,y)$ is an efficient coalition, and coalition $W =V(P;x,\dots,y)$, which contains it, is also efficient!

\

\noi 2. Conversely, for each $i$ belonging to $W$, there exists at least one individual choice, $c_i$, such that $T(c_i,x_i,\dots,y_i)$ is true. We then construct a collective choice, $c$, taking arbitrarily, for example, $c_i = 0$ when $i$ does not belong to $W$  [it is clear that this choice is irrelevant]. Then coalition $W$ is contained in $V (c, x, \dots,  y)$ for this choice of $c$. If, therefore, $W$ were efficient, then $V (c, x, \dots, y)$ would also be, which implies (by the induction hypothesis!) that $T (c, x, .\dots, y)$ is true and this, in turn, implies that the statement \lq\lq exists $z, T (z, x, \dots , y)$" is true which means, precisely, that $P (x, \dots , y)$ is true. \qed

\

\su{Remark}- Observe again the very last part of the previous proof. If unanimity, or simply, an effective coalition is satisfied with a same classical, real object, $e$ for example, to approve the statement $T (e, x, . . . , y)$, there is no hesitation to designate it as a choice of the community. Otherwise, we create the new object $c$, fictitious, \lq\lq nonstandard" that will represent the collective choice. An ideal object.

\

\noi Mathematicians for a long time, handle \lq\lq ideals" in multiple and various forms: Kummer ideals, points at infinity, imaginary numbers ... the list is inexhaustible.

\

\noi From this lemma immediately follows, as a corollary, the next key result on permanence of statements.

\

\su{17 The transfer theorem} {\bf A classical statement is true if and only if the assembly decides that it is true.}

\

\noi Applied to statements of the form (not $P$), this result, assuredly, says that a classical statement is false if and only if the assembly so decides.

\

\noi The echo, returning back, says: \lq\lq We can perfectly well decide by a vote if a theorem is true or false. Just look at this fictitious world, seen through the eyes of members of an infinite assembly!"

\

\su{18 An infinite assembly}

\

\

An infinite assembly, this does not exist, some would say. This can, however, easily be created by a mathematician.

\

\noi Imagine an infinite countable assembly whose members would be nat- ural numbers, $1, 2, \dots , n, \dots$ . Let such an assembly vote according to a non-dictatorial ultrafilter. Such ultrafilters exist, we already said so. The system not being dictatorial, no member alone can be an efficient coalition by himself. More generally, under condition {\bf U2}, no finite coalition can be efficient, so that the opposing (cofinite) coalition is efficient. Ask the assembly to elect an integer. Imagine that each member votes for himself. What can be said about the winner? Denote him by $\omega$. He is the collective choice, the assembly decided he was the winner. He is an integer since the assembly is unanimous to say so. Compare this choice to the integer 7, for example. Only the first seven members of the assembly would say he does not exceed 7. All others, a majority, agree to say that $\omega$ exceeds 7. The same would happen if one compared $\omega$ to any other integer $n$. We would always have an efficient coalition say that $\omega$ exceeds $n$. And the assembly is sovereign! Thus, by decision of the assembly, we can declare that $\omega$ is an \lq\lq integer greater than all natural numbers." He well deserves the name of \lq\lq infinitely large integer". The winner $\omega$  is {\bf an infinite integer}.

\

\noi Similarly, voting for a real number, if each member $n$ of the assembly chooses the fraction $1/n$, the collective choice of the whole assembly (nothing else but $1/\omega$) is a \lq\lq positive, \lq rational' number, smaller than any real positive number". So it is an infinitely small number (not zero), i.e., a {\bf positive  infinitesimal}.

\

\noi Here, we are in the domain of the theory of renovated infinitesimals, alias {\it nonstandard analysis}, of which {\bf A. Robinson} gave, in 1961, perfectly rigorous foundations which nobody expected anymore because people stopped believing they could exist [{\bf ROBINSON, A.}, Non Standard Analysis, North-Holland, Amsterdam, 1966)]. This \lq\lq non-standard" analysis brought to heights at its beginnings, much reviled subsequently did not, probably, deserve all this honor nor this unworthiness! But this, altogether, is another story ...

\

\noi \lq\lq By its choices and successive decisions, the assembly produces new mathematical objects, determines their properties, decides of their mutual relations, states noteworthy theorems which overflow \lq classical' mathematics and, particularly, contains all the elements needed for a sound theory of infinitesimals." [\'Elections, ultrafiltres, infinit\Žsimaux ou le paradoxe de Condorcet, loc. cit.].

\

\noi \lq\lq Through this, and even if at each step we do not reach unanimity, we find unanimity again for statements about conventional objects. This \lq underground' journey of thought is similar to the process which, through imaginary quantities, leads back to real identities (Cardan formulas! Moivre formulas!). One, no doubt, will evoke, too, the per- manence of identities by analytic continuation and monodromy.

\

Passing through fiction does not alter reality in any way when reality reappears." [Condorcet et les ultrafiltres, loc. cit.].

\

\noi Not only such a process alters in no way reality, but it enriches it, greatly.

\

{\bf \lq\lq  1785 \ Condorcet,

 \ \ 1952 \ Guilbaud,

\ \ 1937 \ Cartan,

\ \ 1951 \ Arrow, 

\ \ 1955 \ \L o\'s,

\ \ 1951 \ Bourbaki,

\ \ 1961 \ Robinson.}

\

\noi Contemplating these dates, one starts a dream. A red thread seems to connect them through the years. It goes, it comes, goes back in time, runs and comes back again, and weaves a weft going from Condorcet Paradox to infinitesimals. Can one help imagining that the Marquis, unconsciously or not, have had a premonition of them, that he made, one night, \lq this strange and striking dream'? " [\'Elections, ultrafiltres ... loc. cit.].

\

\su{19 The second degree}

\

\

Complicating, one can imagine several assemblies $I_p$, each equipped with a voting system represented by an ultrafilter $\cal U_p$. Each assembly deliberates and determines itself independently of the others. In order to come out with a joint decision for all these assemblies, it is agreed to resort to a voting system, \lq\lq for a synthesis", defined by an ultrafilter $\cal U$  on the set $P$ of all indices $p$, that is, on the large gathering of the respective spokesmen $p$ of these assemblies. The procedure is as follows: For a given question, the collective response of the family of assemblies is \lq\lq yes" when the set of all those $p$ whose assembly $I_p$  answered \lq\lq yes" is an efficient coalition of the assembly $P$ with respect to ultrafilter $\cal U$.

\

This is a type of a two degrees voting procedure, reminiscent of the elections for the Senate in France.

\

\noi Let us simplify. A moment's reflection shows that one can do without a second degree in this procedure. Just consider the union $I$ of all the assemblies $I_p$ and define the set $S$ of efficient coalitions in the assembly $I$ by the following rule: Coalition $K$ is efficient in $I$ whenever the set $\{p \in P : K \cap I_p \in U_p\}$ belongs to $\cal U$.

\

\noi This one degree system is equivalent to the two degree system. In fact, two conditions are assumed: 1) If an individual $i$ belongs to two or more different assemblies $I_p$, his choices are supposed to always be the same in each of these assemblies. 2) The spokesmen are to comply, of course, with the imperative mandates of the assemblies that delegate them.

\

\noi The ultrafilter $\cal S$ defined on the set $I$ for this voting system is the ultrafiltered sum with respect to $\cal U$ of the family of ultrafilters $\cal U_p$. This notion was introduced, explicitly, by {\bf G. Grimeisen} in 1960 for special purposes in general topology [{\bf GRIMEISEN, G.}, Summation von Filtern und iterierte Grenzprozesse. I, Math. Ann. {\bf 141} (1960) 318-342; and II, ibid. {\bf 144} (1961) 386-417)].

\

\noi Let us iterate. Then, the next idea comes naturally. An assembly $I$ governed by a constitution (or voting system) $\cal U$, creates by its choices (as we have seen) a world of new fictitious objects (mathematical objects ... or others!) Consider a second assembly $J$, governed by an ultrafilter $\cal V$ which takes a look [supervises] this new world created by the assembly $I$, the new world that already contains the real world, of course! The idea is not so crazy. One inevitably thinks of second instance courts in France.

\

\noi As above, a moment of reflection, slightly longer, shows that, here too, a simplification occurs.

\

\noi The overall role of these two assemblies can easily be vested to a single assembly $I \times J$ consisting of all pairs $(i, j)$ with the following rule: A coalition $K$ in $I\times J$ is efficient whenever the coalition $L$ of all members $j$ of $J$ for whom the set $K(j) = \{i \in I : (i,j) \in K\}$ belongs to $\cal U$, belongs itself to $V$, that is to say, whenever we have
$$L:= \{j \in J :\{i\in I :(i,j)\in K\}\in \cal U\}\in \cal V.$$
The ultrafilter $\cal W$ of efficient coalitions thus defined on the product set $I\times J$, is the ordinal product of ultrafilter $\cal U$ by ultrafilter $\cal V$; their roles are clearly not symmetric. This operation, as we know, as we see, reduces to an ultrafiltered sum when \lq\lq splitting" the assembly $I\times J$ into \lq\lq horizontal" sub-assemblies $I_j = \{(i, j) :  i \in I\}$.

\

\noi In other words, we let assembly $I$ vote \lq\lq $J$ times", and every member $i$ can vary his choices, in each vote, however he wishes.

\

\noi The world produced by the superposition of the two assemblies, $I$ and $J$, is the same as the world of the product assembly $I \times J$.

\

\su{20 The second order and the double ultrapower}

\

\

Very early, in the use of infinitesimals, the need was felt for several orders of infinitely small and several orders of infinitely large numbers. Upon introduction of nonstandard analysis, anew, the desire to have those two scales has manifested itself. Several means to achieve it had been suggested. It seems to me that there is no need to go too far nor complicate simple issues. To reach second order, the double ultrapower is enough. To achieve higher orders, one can just iterate, and instead of two assemblies, take as many as needed! Let them then vote successively, \lq\lq hierarchically", as shown above. This will provide all levels of orders desired, with no effort, no changes in the manner the fictitious worlds are looked at [HADDAD, L., La double ultrapuissance, {\it S\Žminaire d'Analyse}, Universit\Ž de Clermont II, ann\Že 1987-1988, ex- pos\Ž $n^\circ$ 24].

\

\su{21 Ultralimits} One more word about this subject to say the following. The utralimit corresponds to the superposition of an infinite sequence of assemblies, \lq\lq potential", not \lq\lq actual". Curiously enough, here, the reduction of the sequence to a single assembly is no more appropriate!

\

\su{22 Condorcet and ultrafilters}

\

\

\lq\lq It is very unlikely to be belied if one thinks that the Marquis de Condorcet was unaware of ultrafilters. Similarly, it is quite likely that Henri Cartan, inventing them, was thinking more \lq to formulate the notion of limit in all generality', rather than avoid the Marquis Paradox.

\

\noi Condorcet and ultrafilters. However, the encounter was inevitable. It owes less to chance than to necessity!

\

\noi The metaphor that results allows any layman, as we have seen, to enter the nonstandard chapel without prior initiation.

\

\noi A metaphor is properly a transportation, a transfer of meaning. This can lead far away ... "[Condorcet et les ultrafiltres, loc. cit.].

\

\su{23 Provisional Conclusion}

\

\

There had been thoughts given to measure \lq\lq the cohesion of a community" in inverse proportion to the distortion between real frequencies of situations \ˆ la Condorcet and theoretical probabilities. To create a typology according to the various configurations of coalitions. The idea of measuring abstractly the cohesion of a community using the frequency of situations \ˆ la Condorcet is attractive. However, I have seen few texts putting it into practice.

\

\noi In human societies, one will probably have to learn how to cope with the Condorcet effect. However, it is curious to hear so seldom, and so discreetly, talks about it. We cannot but feel, at times, that those who know stay careful not to warn those who do not yet know.

\

\noi The problem of general interest, or average opinion, in addition to its political social, economic, psychological, legal, logical, ... and mathematical facet, certainly presents a philosophical aspect. Of course, \lq\lq Shoemaker, not above the sandal" ({\it Sutor, ne supra crepidam}). But one cannot help thinking that Condorcet's Paradox shall always remain a tight knot in the very heart of the pursuit of general interest and common good.

\

\noi As game theory teaches, two completely informed players always have, theoretically, a resting place to retrench. With three players, and on, no possible retreat, plays get animated, games become dangerous.

\

\noi As we know, underground trains of thought often borrow on metaphors and metonymies. The conscious activity too.

\

\noi From a mutual fund, we take our own pictures. To each his favorite. A multitude of looks laid on the same object!

\

\noi If we can associate infinitesimals to the resolution of conflicts, I will, myself, see more than a chance encounter, an enrichment, so minute be it, of our imagination, a small \lq\lq'bunch of pictures" [Denjoy's expression].

\

\su{24 Conditions for a consistent ranking}

\

\

We will not leave the subject of voting systems without addressing a question that seems to have much preoccupied theorists in economy after the work of Arrow: To find {\it conditions} for the aggregation of individual rankings to lead to a collective consistent classification. {\it Conditions} that must be imposed on a set of individual rankings [{\it profiles} as they sometimes are called], not to the voting system itself, of course.

\

\noi Here is the {\bf problem}.

\

\noi We are given an assembly $A$ and a number of candidates $a,b,c,\dots$. Each member of the assembly makes his own ranking, putting these candidates in the order of his preferences, for example, $a > b > c > \dots$. It is assumed that these rankings are linear orders [avoiding ties]. The assembly then proceeds to a comparison of the candidates, two by two, by successive votes. This results into a collective ranking. The issue is to find conditions (on the set of all individual rankings) so that this collective classification be, itself, a linear order, i.e., contains no cycle.
It will be noticed, by the way, that the question needs only be settled in the reduced case where there are three candidates only. Indeed, the collective ranking is a linear order if and only if the induced order on each one of the triples of candidates, $\{x,y,z\}$, is linear. \{A binary relation on a set is a linear order if and only if its restriction to each of the triples of this set is a linear order! This is clear.\}

\

Let us mention the recent text of {\bf Elsholtz Christian} and {\bf Christian List}, A Simple Proof of Sen's Possibility Theorem on Majority Decisions, {\it Elemente der Mathematik}, 60 (2005) $n^\circ$ 2, 45-56. Here is the authors summary.

\

\lq\lq Condorcet's voting paradox shows that pairwise majority voting may lead to cyclical majority preferences. In a famous paper, Sen (1966) [Sen, A.K., A Possibility Theorem on Majority Decisions, Econometrica {\bf 34} (1966) 491-499, reprinted in Sen, A.K. (1982) Choice, Welfare and Measurement, Oxford, Blackwell] identifies a general condition of individual preference orderings, called triplewise value-restriction, which is sufficient for the avoidance of such cycles. This note aims to make Sen's result easily accessible. We provide an elementary proof of Sen's possibility theorem and a simple reformulation of Sen's condition. We discuss how Sen's condition is logically related to a number of precursors. Finally, we state a necessary and sufficient condition for the avoidance of cycles, and suggest that, although there is still some logical space between that condition and Sen's sufficient condition, Sen's condition cannot be further generalized in an appealing way".

\

\noi It seems, according to Elsholtz and List, that Sen's (sufficient) con- dition is the best that has so far been found. Let us emphasize the last sentence: \lq\lq{\it   ... we state a necessary and sufficient condition for the avoidance of cycles, and suggest that, although there is still some logical space between that condition and Sen's sufficient condition, Sen's condition cannot be further generalized in an appealing way}". All this takes place in the domain of the classical majority voting system.

\

\noi We will show how to significantly improve these results, make them perfectly clear and legible, in the more general context of generalized voting systems, and how to obtain a condition, both {\bf necessary} and {\bf sufficient}, quite {\it attractive} and {\it appealing}.

\

\noi  A few words,  first, to recall what Sen's \lq\lq triplewise value-restriction" is. It is the following condition:

\

(SEN) Whatever the triple $\{x,y,z\}$ of candidates, there is a rank $r \in \{1,2,3\}$ such that one of these candidates, $t \in \{x,y,z\}$ is not ranked $r$ by any of the voters.

\

\su{25 A general framework}

\

\

Take any voting system, $(A,\cal E)$, for an assembly $A$, where $\cal E$, the set of efficient coalitions, satisfies conditions {\bf C1} and {\bf C2}]. [Of course, as already said, the classical majority system fulfills these conditions (provided a clause is added, that of a {\it casting vote}, for example).]

\

\noi For any given three candidates, $a, b, c$, only six consistent rankings exist. Label them as follows using the elements of the cyclic group $\Z/6\Z$ of order 6:

\

$1 \ \ a>b>c$ 

$2 \ \ a>c>b$ 

$3 \ \ c>a>b$ 

$4 \ \ c>b>a$ 

$5\ \  b>c>a$ 

$6 \ \ b>a>c$.

\

\noi This labelling has the following peculiarities. Rankings $p$ and $p + 1$ always have either a first same candidate or a last same candidate. To go from $p$ to $p + 1$, one \lq\lq disturbs" as little as possible these rankings, i.e., one simply swaps the ranks of the two first or of the two last candidates. In a sense, rankings $p$ and $p + 1$ are \lq\lq as close as possible". Finally, rankings $p$ and $p + 3$ are \lq\lq opposite". For example, one has

\

$1 \ \ a>b>c$ 

$4 \ \ a<b<c$.

\

\noi Thus, in particular, the same candidate occupies the second rank in two different rankings, $p$ and $q$, if and only if we have $p + q = 3$, i.e., when the two rankings $p$ and $q$ are opposite.

\

\noi That being said, introduce the following notations: Denote $K(p)$ the coalition of those members of the assembly who choose $p$ for a ranking of the three candidates. Similarly, denote $K(p,q)$ the union of coalitions $K(p$) and $K(q)$, and let $K(p,q,r)$ be the union of coalitions $K(p)$, $K(q)$, and $K(r)$.

\

\su{26 A solution}

\

\

Then introduce the following three conditions.

\

(S) There is a $p$ such that coalition $K(p, p+1)$ or coalition $K(p, p+3)$ is empty.

(T) There is a p such that coalition $K (p, p + 1)$ is efficient.

(V) There is a $p$ such that  both the coalitions $K(p, p + 1, p + 2)$ and $K(p+1,p+2,p+3$) are efficient.

\

(1) Condition (S) is the natural generalization of Sen's condition because, in the majority system, condition (S) is none other than condition (SEN) 
above.

\

(2) Condition (T) is a consequence of condition (S).

\

{\bf Indeed}, to see that (S) implies (T), it suffices to notice this. When $K(p, p + 1)$ is empty, the two coalitions $K(p + 2, p + 3)$ and $K(p + 4, p + 5)$ are opposite so that one of them is efficient. Similarly, when $K (p, p + 3)$ is empty, the two coalitions $K (p + 1,p+2)$ and $K(p+4,p+5)$ are opposite and one of them is efficient. \qed

\

(3) Condition (T) implies condition (V).

\

{\bf Indeed}, if  $K(q, q + 1)$ is efficient, both coalitions $K(q - 1,q,q + 1)$ and $K(q,q + 1,q + 2$) which both contain it are  efficient. Then, just take $p = q - 1$. \qed

\

(4) Condition (V) implies that the collective ranking is linear.

\

{\bf Indeed}, suppose that  $K(p, p + 1, p + 2)$ and $K(p+1,p+2,p+3)$ are both efficient. The ranking $p$ has the form $x > y > z$.

\

There are only two possible cases.

\

In the first case, one has

$p \quad \quad \ x>y>z$ 

$p+1 \ \  x>z>y$ 

$p+2 \ \ z>x>y$ 

$p+3 \ \ z>y>x$.

\

\noi The efficient coalition $K(p, p + 1, p + 2)$ thus imposes the collective preference $x > y$ and efficient coalition $K(p+1, p+2, p+3)$ the collective preference $z > y$. Whatever the collective preference between $x$ and $z$, the collective ranking will always be linear! [It will be either $p+1$ or $p+2$.]

\

\noi In the second case, one has

\

$p \quad \quad \ x>y>z$ 

$p+1 \ \ y>x>z$ 

$p+2 \ \ y>z>x$ 

$p+3 \ \ z>y>x$.

 \
 
\noi Here, the two efficient coalitions $K(p, p + 1, p + 2)$ and $K(p + 1,p+2,p+3)$ force, respectivly, collective preferences $y > z$ and $y > x$ so that the collective ranking is still linear. [This will, again, be either p+1 or p+2.] \qed

\

(5) A small surprise finally comes out from the following result.

\

{\bf Condition (V) is not only sufficient, but also necessary in order that the classification of the three candidates $a,b,c$, be consistent.}

\

{\bf Indeed}, suppose the collective ranking of the three candidates is linear. Even if we have to change the names of the candidates, we can assume that this ranking is $a > b > c$. The coalition of voters for whom $a > b$ is efficient: This coalition is none other than $K(1, 2, 3)$. Similarly, the coalition of voters for whom $a > c$ is efficient and that is coalition $K (6, 1, 2)!$ So that condition (V) is satisfied. \qed

\

It would have been noticed that condition (V), in essence, says that individual choices must be neighboring and not too much \lq\lq scattered" which was predictable.

\

(6) The necessary and sufficient condition (V) and, even more, the sufficient condition (T), certainly have some character of simplicity. Of course, one can easily translate each of these conditions, for the case of the conventional majority voting system, in terms of voting numbers in different coalitions. It must then be limited to the case where voters are odd in number or accept the complications due to one of the additional terms like \lq\lq casting vote". {In the case of a majority system for an odd number of voters, for example, condition (T) is written just like this: There is a $p$ for which we have $|K(p, p + 1)| > |K(p, p + 1)^c|$.}

\

(7) We have here dealt with the cases where rankings do not include ties. The general case where rankings are preorders can have a similar treatment, slightly more complex, of course.

\

\su{27 Remarks}

\

\

(1) The issue arises as to the precise place condition (T) has between the two conditions (S) and (V). Simple examples, in the majority system itself, can show that (T) is neither equivalent to (S) nor to (V).

\

In order to see that condition (T) does not imply (S), one has just to look at the following situation: 5 voters out of whom 3 choose ranking 1, the others two choosing, respectively, ranking 3 and 5, which means  that $|K(1)| = 3,|K(3)| = |K(5)| = 1$ and $K(2) = K(4) = K(6) = \vide$. Thus, none of coalitions $K(p, p + 1)$ and $K(p, p + 3)$ is empty while coalition $K (1, 2)$, for example, is a majority. Condition (T) is satisfied not (S). \qed

\

To see that condition (V) does not imply (T), look at the following situation: 5 voters choose, respectively, the first 5 rankings, $1,2,3,4,5$, no one chooses $6$, i.e., each of coalitions $K(1), K(2), K(3), K(4)$, and $K(5)$ is a singleton while $K(6)$ is empty. Thus, each of the coalitions $K (p, p + 1)$ contains at most 2 members so that none is a majority. However, coalitions $K(1,2,3)$ and $K(2,3,4)$, each, contains 3 members and, therefore, both are majorities. This means that condition (V) is satisfied but not (T). \qed

\

(2) Since (T) implies (V) which in turn implies lack of inconsistency, this condition (T), a fortiori, implies the consistency of the collective ranking.
Here is a direct proof which thus also proves, in a very simple way, that \lq\lq Sen'"s condition (S) is sufficient to ensure consistency of the collective choice.

\

{\bf Indeed}, if a coalition $K (p, p + 1)$ were efficient, it would dictate the choice of the first or of the last candidate in the two rankings, $p$ and $p + 1$, depending on the cases, and the ranking of the two other candidates would not introduce any inconsistency. [For example, if $K(1,2)$ were efficient, collective preference would be $a>b$ and $a>c$. It little matters, then, if we had $b > c$ or $c > b$ , the collective ranking will always be linear.] \qed

\

\su{28 The ultrafilter: A conciliator} Out of the cacophony that would result if each member of the infinite assembly $1,2,3,...,n,...$ gives his own opinion, the utrafilter $\cal U$ derives a unique sound that harmonizes all these opinions, so diverse and varied.

\

\noi It can also have a similar role in classical analysis, as we shall readily see. We also take time for some reminders.

\

\su{29 The ultrafilter: A moderator}

\

\

The behavior of a real sequence which is not convergent can be quite erratic, very chaotic, as we know.

\

\noi Let $\cal U$ be a nontrivial ultrafilter on the set $\N = \{0,1,2,...,n,...\}$ of natural numbers. Each real sequence $x = (x_n)$ has a limit, $\lim_{\cal U} x$, with respect to this ultrafilter. \{Recall the following: The limit of a real sequence $x = (x_n)_{n\in \N}$ with respect to an ultrafilter $\cal U$ on $\N$ is the sole element $r$ of the completion $\overline \R = \R \cup \{\pm\infty \}$ such that for each neighbourhood $V$ of $r$, there exists a member $X\in \cal U$ such that $\{x_n : n\in X\}\inc V$. The concept of limit with respect to an ultrafilter is classical. It is to be found in every textbook on topology.\} This limit is a real number when the restriction of the sequence $x$ to a $K \in  \cal U$ is bounded. Otherwise, it is equal to $+\infty$ or $-\infty$, depending on cases. Imagine the assembly $\N$ equipped with the voting system defined by the ultrafilter $\cal U$. Let $e$ be the winner relative to a given collective choice $x = (x_n)$. From the assembly's point of view, the winner $e$ is \lq\lq a real number", let us say {\it hyperreal} to avoid confusion. When the limit $\lim_{\cal U} x$ is finite and equal to $r$, the elect, $e$,  the hyperreal number,  is {\bf infinitely close} to the real number $r$ in the (very precise) sense that the difference $e - r$ is {\bf infinitesimal}. Otherwise, the hyperreal $e$ is infinitely large, either positive or negative, depending on cases. Thus, in particular, every finite hyperreal is infinitely close to a real number which is its {\it shadow}
in the real world, somehow. The shadow of an infinitely large hyperreal is, of course by definition, either $+\infty$ or $-\infty$, according if it is positive or negative.

\

[In this context, and as far as I know, it was {\bf H. Jerome Keisler} who introduced the word {\it hyperreal}. As for the word {\it shadow} [ombre], so appropriate, I know that we owe it to {\bf Georges Reeb}.]

\

\noi The ultrafilter {\it moderates} the behavior of sequences and to each object it creates in the fictitious world [\lq\lq hyperreal"] attaches an infinitely close object in the real world [its \lq\lq shadow"].

\

\noi Notice also the following. The function $f_{\cal U}$ defined as $f_{\cal U} (x) = \lim_{\cal U} x$ is a positive linear form on the space $B(\N)$ of bounded real sequences.

\

\su{30 Banach (generalized) limits}

\

\

As we well remember, {\bf Stefan Banach} introduced the notion of a generalized limit in his famous book: Th\Žorie des op\Žrations lin\Žaires, Warszawa, 1932, a few years before the introduction of the concept of ultrafilter. He expresses himself, in particular, as follows (on p.34):

\

\lq\lq \`A toute suite born\Že \{$\xi_n$\} on peut faire correspondre un nombre $\mathrm{Lim}_{n\to\infty}\xi_n$ de fa\c con que les conditions suivantes ... soient remplies."

\

[ To any bounded sequence \{$\xi_n$\} can be attached a number $\mathrm{Lim}_{n\to\infty}\xi_n$ such that the following conditions ... hold.]

\

\noi The conditions which the operation Lim, defined on the space $B(\N)$, must satisfy (expressed in today's language) are: Lim is a positive linear form, such that $\mathrm{Lim}_{n\to\infty}\xi_{n+1}= \mathrm{Lim}_{n\to\infty}\xi_n$ and $\mathrm{Lim}_{n\to\infty}1=1$.

\

\noi \{The Lim operations which satisfy these conditions were called \lq\lq Banach generalized limits."\}

\

\noi Of course, Banach establishes the existence of these operations with the help of his theorem on \lq\lq the extension of the linear functionals" bounded by a \lq\lq semi-norm", as we would say nowadays.

\

\noi Previously, {\bf Ren\Ž Baire} had introduced the concepts of upper and lower limits of a real sequence which satisfy the following inequalities:
$$\underset n \inf(x_n) \leqs\underset{n\to \infty}   \liminf \ x_n \leqs \underset {n\to \infty} \limsup \ x_n \leqs  \underset n \sup(x_n).$$

\

\noi Of course, any positive linear form $f$ defined on the space $B(\N)$ satisfies, for each real sequence $x = (x_n)$, the following inequalities
$$\underset{n\to\infty} \liminf \ x_n \leqs  f(x) \leqs \underset{n\to\infty}  \limsup \ x_n.$$

\

\noi Therefore $f$ is, in particular, continuous for the norm $||x|| = \sup_n |x_n|$, as we know.

\

\noi The condition $\mathrm{Lim}_{n\to\infty}\xi_{n+1}= \mathrm{Lim}_{n\to\infty}\xi_n$ says that the operation Lim is invariant under shifts of sequences, which amounts to say that any finite number of terms of the sequence can be removed without changing the value of its generalized limit.

\

\noi As for the {\it normalization} requirement, $\mathrm{Lim}_{n\to\infty}1=1$, it is not really essential, but to exclude the null form. In fact, given any positive linear non zero form $f$ defined on the space $B(\N)$, invariant under shifts, the function $g(x) = f(x)/f(1)$ is a generalized Banach limit.

\

\noi Let us take one step more. Let $f$ be any positive linear form, defined on the space $B(\N)$, and such that $f(1) = 1$. For each real bounded sequence $x = (x_n)$, first set $t_n = (1/n) \sum_{1\leqs k\leqs n} \ x_k$ then $t = (t_n)$ and $g(x) = f(t)$. The function $g$ is then a generalized Banach limit. \{It is a straightforward exercise to check that $g$ has the the required properties.\}

\

\noi In other words, the emphasis is on the set of all positive linear forms on the space $B(\N)$.

\

\noi That being said, let us go back to a nontrivial ultrafilter $\cal U$ on $\N$ and, for each sequence $x \in B(\N)$, set $t_n = (1/n) \sum_{1\leqs k\leqs  n} \  x_k$  then $t = (t_n)$, and finally $g(x) = f_\cal U(t) = \lim_\cal U(t)$. The function $g$ is a generalized Banach limit!

\

\su{31 Ultrafilter spaces}

\

\

Let $E$ be any set and denote  $\Upsilon(E)$ the set of all ultrafilters on $E$. [This should be read \lq\lq upsilon" of E.] We already know, indeed, that
each point $x\in E$ defines a trivial ultrafilter $\cal U_x$. We thus have an injective map $x \mapsto \cal U_x$ from the set $E$ into $\Upsilon(E)$. This injection is used to identify the set $E$ with the subset of trivial ultrafilters in $\Upsilon (E)$. This identification, point to point, is safe, except in exceptional cases.

\

(1) If the set $E$ is finite, all ultrafilters on $E$ are trivial so that we have $\Upsilon (E) = E$.

\

(2) If the set $E$ is infinite and if we use a set theory with axiom of choice, we already know that nontrivial ultrafilters exist on $E$ and we can, in fact, prove that there are {\it quite a lot} of them. To be specific, if $\mathrm{card}(E) = \alpha$ is infinite, then $\mathrm{card}(\Upsilon(E)) = 2^{2^\alpha}$.

\

(3) It is also known that there are models of set theories in which infinite sets exist but nontrivial ultrafilters do not (assuredly, theories without the axiom of choice). In other words, in these models, we always have $\Upsilon (E) = E$, whatever the set $E$ is.

\

(4) To each subset $F \inc E$ is also associated the set $\Upsilon(F)$ of ultrafilters on $F$. Each ultrafilter $\cal V$ on $F$ is the trace of a unique given ultrafilter $\cal U$ on $E$. Identifying $\cal V$ to $\cal U$ , we can also, safely, identify $\Upsilon (F)$ to the corresponding subset of $\Upsilon(E)$.

\

(5) That being said, the set $\{\Upsilon (F) : F \inc E\}$ is a basis for the open sets of a topology on the set $\Upsilon (E)$. This {\bf canonical} topology turns the set $\Upsilon(E)$ into a compact space, the {\bf Stone space} of ultrafilters on $E$. This space is very usually denoted $\beta E$. It is, as one can prove, the {\bf Stone-\v Cech  compactification} of $E$ considered as a space with its discrete topology. \{There are many ways to obtain compactifications. The method that consists in using ultrafilters is one of the most natural. It is well known. We will see, below, how $\beta E$ is the spectrum of the ring $\cal P(E)$.\}

\

(6) Let us go a little bit backwards [see paragraph {\bf 19}] to have a new look at the family of ultrafilters $(\cal U_p)_{p\in P }$ on a set $I$ and the ultrafilter $\cal U$ on the set $P$ of indices. The [Grimeisen] ultrafiltered sum relative to $\cal U$ of the family $(\cal U_p)_{p\in P}$ simply appears as the limit, in the compact Stone space $\beta I$, of the family of ultrafilters relative to the ultrafilter $\cal U$.

\

(7) Spaces of ultrafilters play a considerable role in general topology. They hold an important place in most of the treaties. [Bourbaki, however, relegated them as exercises which form, one must confess, a non-negligible part of his treaty.] It is known, for example, that all compactifications of a completely regular space $E$ are quotients of the compact Stone space $\beta E$.

\

(8) The Stone space $\beta E$ has also another interpretation. It is the {\bf spectrum} of the Boolean ring $\cal P(E)$ of the subsets of $E$. Indeed, the prime ideals of this Boolean ring, which are exactly its maximal ideals, are in a bijective correspondence with ultrafilters on E. More precisely, let $\cal M$ be a maximal ideal in the Boolean ring of subsets of $E$. The set
$$\cal U = \{K : K \inc E \ \mathrm{and} \ (E \stm K) \in \cal M\}$$
is an ultrafilter on $E$ and, of course, we also have the {\it dual}
relation
$$\cal M = \{K : K \inc  E \ \mathrm{and}\ (E \stm K) \in\cal U\}.$$
All those things, and many more, are classical and well known.

\

\noi Let us say a few words of a somewhat disregarded use of ultrafilter spaces and which, perhaps, it is interesting to be aware of, all the same.

\

\su{32 Nasses [creels]}

\

\

Le$E$ be any topological space whose set of open subsets is $\cal O$. Define a graph $T$ in the product space $\Upsilon(E) \times \Upsilon (E)$ as follows:
$$(\cal U,\cal V)\in T \iff \cal U\cap \cal O\inc V.$$

The graph $T$ was named the {\bf  nasse} of the topological space $E$. [{\bf HADDAD, L}. Une repr\Žsentation des topologies, C. R. Acad. Sci, Paris, 255 (1962) 2702-2704, and Sur quelques points de topologie g\Žn\Žrale. Th\Žorie des nasses et des tramails, {\it Ann. Fac. Sci. Univ. Clermont}, $n^\circ$ {\bf  44} (1970) fasc.7, 3-80.]

\

\noi The graph $T$ thus is a binary relation on the set of ultrafilters on $E$ and, as usual, in order to say that we have $(\cal U,\cal V) \in T$ , we can write either one of the following:
$$\cal U \ T\ \cal V\  \mathrm{or} \  \cal V \overset{-1} T\cal U \  \mathrm{or}\ \cal V\in T(\cal U) \  \mathrm{or}\ \cal U\overset{-1}T(\cal V)$$
where $\overset{-1}T$ is the opposite, (reciprocal or inverse) relation.

\

\noi Without any detour, one verifies that $T$ is a reflexive and transitive relation, i.e., we always have
$$\cal U \ T \ \cal U$$
$$\cal U \ T \ \cal V \ \mathrm{and} \ \cal V \ T \ \cal W \implies \cal U \ T \ \cal W.$$
\{Simply, starting with the definitions and taking the direct path, one sees that relation $T$ is reflexive since $\cal U \cap \cal O \inc  \cal U$. Moreover, $\cal U \cap \cal O \inc \cal V$ and $\cal V\cap O\inc \cal W$ clearly imply $\cal U\cap \cal O\inc \cal W$ whence $T$ is transitive.\}

\

The nasse $T$ is a {\bf preorder} on the set $\Upsilon(E)$ of ultrafilters. It contains all the information about the topology of the space $E$ and thus reduces, in a sense, the study of topologies to that of preorders.

\

\noi When the set of all ultrafilters, $\Upsilon(E)$, is equipped with its canonical topology, i.e., when one uses the Stone space $\beta E$, one perceives that the nasse $T$ is a {\it closed} subset in the product space $\Upsilon (E) \times \Upsilon (E)$. Moreover, it has the following peculiarity: The preimage by $T$ of each open subset of $\beta E$ is again an open subset of $\beta E$. 

\noi That is,
$$\cal A \ \text{is an open subset of} \ \beta E \implies \overset{-1}T(\cal A) \ \text{is an open subset of} \ \beta E,$$
as in the case of continuous mappings. Those binary relations which thus mimic continuous mappings, {\bf Choquet} called them {\bf mi-open} relations [relations {\bf mi-ouvertes}]. The nasse $T$ is thus {\bf a closed and mi-open preorder on the Stone space $\beta E$}.

\

\noi What is remarkable about those conditions is that they characterize the nasses of topologies, in the sense that each {\bf closed and mi-open preorder on the Stone space $\beta E$} is the nasse of a topology on the set $E$. There is a bijective correspondence between the set of topologies on $E$ and the set of closed and mi-open preorders on $\beta E$.

\

\noi More generally, a nasse is defined to be a reflexive binary closed relation on the space $\beta E$. A binary relation $R$ is said to be {\it idempotent} when $RR = \overset 2 R = R$. So nasses of topologies are, precisely, the {\bf idempotent and mi-open nasses}.

\

\noi \lq\lq Remove the idempotence and fall essentially on Choquet pretopologies. Take out the mi-open condition and find the topogenous orders of Cs\'asz\'ar.
Add the condition of symmetry and obtain a useful representation of proximities." [See \lq\lq Quelques points de topologie g\Žn\Žrale. Th\Žorie de nasses et des tramails", loc. cit..]

\

\noi Nasses cover a wide-range of classic topological structures and their best known generalizations. In the text quoted above, one will find a complete (enough) treatment of these questions. To give a small foretaste, we will add this.

\

\noi Given a topological space $E$ and its nasse $T$, the necessary and sufficient condition for the space to be {\bf normal} is that we have $T \overset{-1}T \inc \overset{-1} T T$ .
Similarly, the space $E$ is {\bf extremally disconnected} if and only if we have $\overset{_1} T T \inc T \overset{-1}T$.

\

\noi Let us open one last parenthesis in this paragraph to say the following. Applied to the case of finite sets, the introduction of nasses shows, in particular, that topologies on a finite set are in a bijective correspondence with preorders. This fact, quite simple, seems to have been {\it rediscovered} again and again, researchers, encountering it for the first time, could think, legitimately, be the first to have observed it. We can say in a sense that the theory of nasses is an extension of this simple fact to more complex topological structures.

\

\su{33 Choquet and ultrafilters}

\

\

Before concluding, I still like to mention, specially, the quite remarkable use Choquet did of ultrafilters in his work. Alas! I will not have the time nor [above all] the sufficient skills to give a satisfactory overview. Also, I shall content myself with an overview that will focus on only one of the many highlights of all the impressive uses he made of them. What I am going to talk about is a tiny part of a very dense, and essential, text: {\bf  CHOQUET G}., Convergences, {\it Ann. Univ. Grenoble, Sect. Sci. Math. Phys}., {\bf 23} (1947-1948) 57-112.

\

\noi As we know, {\bf Gustave Choquet} is {\bf  Arnaud Denjoy}'s student who, himself, wrote down the superb lessons that {\bf  Ren\Ž Baire} gave at the Coll\ge de France [Le\c cons sur les fonctions discontinues, Gauthier-Villars, Paris, 1905]. Among the three, the filiation is clear and, some day, a historian of mathematics will surely bother to extricate the themes and essential episodes of this filiation then secure them through writing, so the oral tradition will live on.

\

\noi Take any topological space $E$. Take a family $\cal X = (X_i)_{i\in I}$ of subsets of $E$ indexed by a set $I$, and a filter $\cal F$ on $I$. Introduce the concepts of upper and lower limits of the family $\cal X$ relative to the filter $\cal F$ as follows. [As we see, inspiration comes from afar.]

\

\noi Begin associating to the filter $\cal F$ its {\bf grill}
$$\cal G=\{Y :Y \inc E, \  X\cap Y \neq \vide \ \text{forall} \  X\in \cal F\}.$$
Denote $\Upsilon < \cal F >$ the set of all those ultrafilters on $I$ that contain [i.e., which are finer than] $\cal F$. Of course, the filter $\cal F$ is the intersection of all those ultrafilters which implies, as can easily be seen, that  $\cal G$, the grill,  is the union of these ultrafilters, $\cal G = \bigcup_{\cal U\in \Upsilon<\cal F >} \cal U$. It can be noticed that the grill $\cal G$  is equal to the filter $\cal F$ if and only if $\cal F$ is an ultrafilter!

\

\noi For each subset $J \inc  I$, set $X_J = \bigcup_{i\in J} X_i$. Let $\overline{X_J}$ denote the closure of $X_J$ and set
$$\underset{\cal F} \liminf \cal X = \underset{J\in \cal G}\bigcap \overline{X_J}$$
$$\underset{\cal F} \limsup \cal X = \underset{J\in \cal F}\bigcap \overline{X_J}.$$

\

(1) As should be, the upper limit is larger than the lower limit:
$$\underset{\cal F} \liminf \cal X \inc \underset{\cal F} \limsup \cal X,$$
they both are equal if $\cal F$ is an ultrafilter, but, sure enough, this is not the only case when they are equal. When both limits are equal, denote $\lim_\cal F \cal X$ their common value and say that the family $\cal X$ has {\bf a limit relative to the filter $\cal F$}. In particular, the family $(X_i)_{i\in I}$ always has a limit relative to each of the ultrafilters on $I$, whichever they are. Notice that, more generally, for each of the filters $\cal F\inc \cal  F'$ on $I$, we have:

$$\underset{\cal F} \liminf \cal X \inc \underset{\cal F'} \liminf \cal X \inc \underset{\cal F'} \limsup \cal X  \inc \underset{\cal F} \limsup \cal X.$$

(2) By their nature, the upper and lower limits are closed subsets of the space $E$ and their values do not change on substituting to the family $\cal X$ of subsets $X_i$ the family $\overline {\cal X } = (\overline {X_i})_{i\in I}$ of the closures of the $X_i$'s. So, attention can be restricted to families of closed subsets of $E$.

\

(3) Denote $\Phi(E)$ the set of all closed subsets of the space $E$. We thus have a notion of limit in the set $\Phi(E)$. To be more specific, let $\cal X = (X_X)_{X\in \Phi(E)}$ be an indexed family whose set of indices is $I = \Phi(E)$ and such that each closed set $X$ is its own index. Given any filter $\cal F$ on the set $\Phi(E)$, its upper and lower limits are defined to be the corresponding limits of the indexed family, and we simply set:
$$\liminf \cal F = \underset{\cal F}\liminf \ \cal X$$
$$\limsup \cal F = \underset{\cal F} \limsup \ \cal X,$$
and, if both limits are equal,
$$\lim \cal F = \underset{\cal F}\lim \ \cal X;$$
Does this turn $\Phi(E)$ into a topological space? That is, can this limit be redefined as a limit in a suitable topology on $\Phi(E)$?

\

(4) In particular, each ultrafilter $\cal U$ on $\Phi(E)$ has a limit in $\Phi(E)$. In order that the limits of filters on $\Phi(E)$ thus defined originate from a topology on $\Phi(E)$, they have, first, to satisfy the following condition:

$$\lim \cal F=F \iff \lim \cal U=F \ \text{for each ultrafilter} \ \cal U\supset \cal F.$$
\

\noi That is, those limits must define a {\bf pseudo-topology} in the sense of Choquet on the set $\Phi(E)$. We prove that this indeed is the case showing, more generally, what follows.

\

(5) Take any family $\cal X = (X_i)_{i\in I}$ of subsets of $E$, indexed by a set of indices $I$, and a filter $\cal F$ on $I$. We then have [see Convergences, loc. cit., page 64]:
$$\underset{\cal F}  \liminf \cal X= \underset{\cal U \in\Upsilon <\cal F >} \bigcap \underset {\cal U} \lim \cal X,$$
$$\underset{\cal F}  \limsup\cal X= \underset{\cal U \in\Upsilon <\cal F >} \bigcup \underset {\cal U} \lim \cal X.$$
Since the grill $\cal G$ is the union $\bigcup_{\cal U \in \Upsilon<\cal F>}\cal U$, the first formula is seen to hold, with little detour. As for the second, in order to show that it also holds, it may be convenient to observe the following result which can be of interest by itself.

\

\noi Given a family $\cal X = (X_i)_{i\in I}$, for each point $x \in E$ and each subset $V \inc  E$, set
$$I(x,V)=\{i : i\in I \ , V \cap X_i \neq \vide\},$$
then
$$\cal I(x) = \{I(x, V ) : V \ \text{neighbourhood of x in E}\}.$$
It is, easily, seen that we have:
$$I (x, U \cap V ) \inc I (x, U ) \cap I (x, V ) \ \text{for all subsets U and V of E},$$
(which proves that $\cal I(x)$ is a filter base on $I$ provided that $\vide$ does not belong to $\cal I(x))$. We then prove, a little less easily, the following two characterizations:
$$x\in \underset{\cal F}\liminf \ \cal X \iff \cal I(x)\inc \cal F,$$
$$x\in \underset{\cal F} \limsup \ \cal X \iff \cal I(x)\inc \cal G.$$
\{The following duality binds a filter $\cal F$ to its grill $\cal G$: 

\noi $(X \in \cal  G) \iff  (X \ \text{meets each} \ Y \in \cal F)$; and $(X \in \cal F) \iff$

\noi $(X \ \text{meets each}\  Y \in \cal G)$. Starting from the definitions, proceed by equivalence, as follows: 

\noi $(x \in \liminf_\cal F \cal X ) \iff (\text{for each} \  J \in \cal  G \ \text{and each neighbourhood V},$ 

\noi $\text{ of x we have} \  J \cap I(x, V ) \neq \vide) \iff (\text{for each neighbourhood V} $

\noi $\text{of x, we have} \  I(x, V ) \in \cal  F) \iff (\cal I(x) \inc \cal F)$, which establishes the first characterization. Do the same for the second.\}

\

That should be enough to prove the second formula.

\
 
\{It is based on the following classical fact: A filter is contained in the grill $\cal G$ of $\cal F$ if and only if it is contained in one of the ultrafilters $\cal U \in \Upsilon < \cal F >$. We thus have $(x \in \limsup_\cal F \ \cal X ) \iff (\cal I(x)\inc  \cal G) \iff (\text{exists} \ \cal U\in \Upsilon <\cal F>  \ \text{such that} \  \cal I(x)\inc \cal U) \iff (x\in \bigcup_{\cal U \in \Upsilon <\cal F>} \lim_\cal U \cal X)$.\}

\

(6) The pseudo-topology thus defined on the set $\Phi(E)$ is not always a topology, not even a  pretopology, as examples can show [see Convergences, loc. cit., page 87]. When $E$ is a Hausdorff space, an important special case, singletons $\{x\}$ are closed sets in $E$, so that $E$ can be identified to a subset of $\Phi(E)$ through an injective map, and, moreover, the pseudo-topology induced on $E$ coincides with the original topology given on $E$. Still more, we then have the following noteworthy results.

\

(7) If the space $E$ is Hausdorff, the pseudo-topology on $\Phi(E)$ is a pretopology if and only if $E$ is {\bf locally compact}. When that is the case, this pretopology is a topology, itself. The subspace $\Phi_0(E) = \Phi(E) \stm \{\vide\}$ is locally compact and $\Phi(E)$ is its Alexandroff compactification. In particular, if $E$ is compact, the space $\Phi_0(E)$ is compact, it is the space of {\it non-empty} compact subsets of $E$.

\

(8) Once wrought, Choquet uses this tool in the study of many-valued relations between two topological spaces [l'\Žtude des relations multivoques entre deux espaces topologiques] which he concludes with a masterful study of relationships between convergence and local uniform convergence, more precisely, between the abstract contingents and paratingents leading to simple statements which are \lq\lq a generalization and geometrization of the results of Baire and those theorems which M. Denjoy put at the base of the theory of functions of a real variable" [loc. cit.].

\

(9) Without going into details, just in order to give a little bit of the flavour of these results, we still add the following.

\

\noi Let $U$ be a metric space, $\Delta$ a compact metric space, $P \inc  E \inc U$ two subsets, and $\delta : (P \times E) \stm (P \times P) \to \Delta$ a map, continuous with respect to its first argument (for each given value of the second). At each point $x \in P$, are defined the contingent $c(x$) of $E$ in $x$ as well as the paratingent $p(x)$ of $E$ in $x$ relative to $P$, associated to the map $\delta$ [that is where it would take much time to enter into details]. Choquet proves the following result.

\

\noi{\bf Theorem.- For each point $x$ in $P$ , except at the points $x$ in a $F_\sigma$  meager subset of $P$, the paratingent $p$ varies continuously and $c(x) = p(x)$.}

\

\noi This is the [utmost] generalization of Baire's classical result about the set of continuity points of pointwise limits of sequences of continuous functions on the interval $[0,1]$!

\

(10) Needless to recall that spaces of closed subsets, $\Phi(E)$, and their subspaces, abound, specially in analysis and geometry. Each is equipped with its canonical pseudo-topology, often thinner than the associated topology and therefore more {\it stingy} on limits. There are spaces of curves, spaces of manifolds, spaces of plane compacta ... There is the subspace $\kappa(E)$ [read \lq\lq kappa" E] of compact susbsets of a Hausdorff space $E$, a subspace of $\Phi(E)$. Let us also recall this. Among these spaces of closed subsets, we must reckon also the many function spaces: The graph of a continuous function on a topological space $E$ with values in a Hausdorff space $F$ belongs to the space $\Phi(E \times F )$ since it is closed in the product space $E \times F$.

\

\su{34 As a fermata}

\

\

Let us get back to our assemblies. Let I be an assembly whose voting system is governed by an ultrafilter $\cal U$. Take a family $\cal X = (X_i)_{i\in I}$  of subsets of a topological space $E$. The assembly transforms $\cal X$ into a collective choice whose representation is an object $M$ in the fictitious world. From its own point of view, the assembly considers $M$ as \lq\lq a subset of $E$", \lq\lq a closed subset" if the $X_i$'s are all closed. The shadow of $M$ in the real world is nothing else but  $\lim_\cal U  \cal X$, in a very precise sense. {\bf But that's another story!}

\

\head{Epilogue}

\

\head{\S1 The paradox of the Marquis de Condorcet}

\

We will try to take a new look at the first two examples of Condorcet.

\su{1 On the first example of Condorcet} Imagine an election where 60 voters must choose one of three candidates, $A, B, C$. A first round gives the following results:

\

$A$ 23 votes

$B$ 19 votes

$C$ 18 votes.

\

If the election is {\it \ˆ l'anglaise} [i.e., a single round, and plurality], then candidate $A$ having the plurality of votes, is elected.

\

If the election is {\it \ˆ la fran\c caise} [i.e., the way the French president is elected], only the two candidates $A$ and $B$, coming top, compete in a second round. So, the supporters of candidate $C$, scattering their votes, decide for the election. [Of course, it is understood that voters do not change their minds and do not abstain in the second round]. Assume that the supporters of $C$ split as follows

\

18 $C >
\bc B > A  \ \ 16\\
A > B \ \ \ 2,
\ec$

\noi that is, out of the 18 supporters of $C$, 16 prefer $B$ to $A$, the other 2 prefer $A$ to $B$. So that $B$ is elected by 19 + 16 = 35 votes /60. 

\

Finally, it can be decided, {\bf \ˆ la Condorcet}, to compare the candidates
pairwise. Let the overall situation be sketched as follows:

\

23 \ $A > 
\bc C>B \ \ 23\\
 B>C \ \ \ 0\\
 \ec
 $

19 \ $B>
\bc C>A \ \ 19 \\

A>C \ \ \ 0
\ec
$

18 \ $C>
\bc B>A \ \ 16\\

A>B \ \ \ 2.

\ec
$

\

\noi Comparing the candidates pairwise, we get the following results: 

\

$B>A$ \ by19+16=35votes/60,

\

$C>B$ \ by18+23=41votes/60,

\

$C>A$ \ by18+19=37votes/60.

\

So, the voters clearly say rather $B$ than $A$, rather $C$ than $B$, and rather $C$ than $A$, each choice being decided by a large majority, which gives the following order $C > B > A$ so that $C$ must be the winner.

\

This example shows that, according to the voting system chosen, the elect is $A, B$, or $C$, respectively. This indicates a significant {\bf sensitivity} of the results to the voting system.

\

But, there is still more, as shown in the following example of Condorcet.

\

\su{2 The second example of Condorcet} Imagine again an election where 60 voters must choose one of three candidates, $A, B, C$. The situation can be sketched as follows:

\

23 \ $A>
\bc B>C \ \  23\\ 
C>B\ \ \  0
\ec
$

\

\

19 \ $B>
\bc C>A \ \ 17\\ 

A>C \ \ \ 2
\ec
$

\

\

18 \ $C>
\bc B>A \ \ \ \ 8\\

 A>B \ \ \ 10.
 \ec
 $

\

\noi Comparing the candidates pairwise, we get the following results: 

\

$A>B$ \ by 23+10=33votes/60,

\

$B>C$ \ by 19+23=42votes/60,

\

$C>A$ \ by  18+17=35votes/60,

\

\noi which is a circular [inconsistent] classification $A > B > C > A$.

\

\su{3 Comments}

\

\

(1) Let us insist again on the fact that the classification thus obtained is circular, inconsistent.

\

(2) {\bf Alliances and coalitions.} In the election {\it  \ˆ la fran\c caise}, $A$ and $B$ are selected at the end of the first round. On the second round, $A$ is elected by 33 votes/60.

\

 We can presume that each voter is well aware of the situation, he knows the marked preferences of all others [as in a kind of game with complete information] and is able to analyze the consequences of all votes, in the light of \lq\lq the rule of the game". Here is what could happen.

\

Among the 19 supporters of $B$, there are $17$ for whom $C$ is a better candidate than $A$. They would be tempted to vote for $C$, in the first round in order he be elected at the outset by 18 + 17 = 35 votes/60.

\

Knowing that, the 23 supporters of $A$ who all prefer $B$ to $C$ would be tempted to offer to the supporters of $B$ an alliance in order to see him elected right out at the first round by 19 + 23 = 42 votes/60.

\

But, among the 18 supporters of $C$, there are 10 that still prefer $A$ to $B$ and who might consider forming an alliance with the supporters of A and make him win at the first round by 23 + 10 = 33 votes/60.

\

Seeing this, the 17 supporters of $B$ who still prefer $C$ to $A$ would offer to add their votes to those of the 18 supporters of $C$ letting him win, so doing, at the first round by 18 + 17 = 35 votes/60. This round dance, has no reason to stop, it could go on for a long time!

\

All this is possible, because there is a majority coalition of 35 voters who prefer $C$ to $A$ and can elect him at the first round. Similarly, There is a majority coalition of 33 voters who prefer $A$ to $B$ and can elect him at the first round. Finally, a majority coalition of 42 who prefer $B$ to $C$, ready to elect $B$ at the first round.

\

How does this round dance stop in real life, in practice? With a lack of any additional data, no one can anticipate nor predict the behavior of the 60 voters. This is probably no surprise to those who read the book by {\bf John von Neumann}  and {\bf  Oskar Morgenstern}, {\it Theory of Games and Economic behavior}, Princeton University Press, Princeton, 1953.

\

(3) I developped this example, one day, in front of young students, future philosophers, who attended my course on {\bf Initiation aux math\Ž\-matiques} [Introduction to Mathematics] in their first University year. Wanting to see their reactions, I asked them what they would do themselves in such circumstances. I was only moderately surprised to hear many answer that they would vote for their first choice at the first round, no matter what could happen, and that any other behavior would amount to \lq\lq skulduggery". I tried, somewhat vainly, to make them feel the difference between {\it skulduggery} and {\it compromise}, showing, among other things, to potential supporters of $B$ that it would be a pity not to improve their terms allying with supporters of $C$ who would still like to see better $B$ elected than $A$! ... Lost efforts ...

\

\noi I had already encountered this form of misunderstanding among many adults discovering for the first time the Condorcet paradox. A lot of \lq\lq stupor" followed by \lq\lq disbelief" leading to sterile \lq\lq denial" ...

\

 Wanting to go further in order to {\it get to the bottom} of the matter, I submitted another example, more striking, so to speak, to my students.

\

\su{4 The Council of Elders} In a remote land, undetermined, and in old times, the life of this City was governed by a Council of Elders who held in their hands the three powers: Executive, legislative and judicial. They established the following law by a two-thirds majority:

\

$p$  Each murderer will be punished.

\

\noi [You can hear: Murderers will be executed, but this time is over.]

\

\noi Some time later they had to judge Untel suspected of being a murderer. Sitting as a court, they listened carefully to the arguments of the prosecutor then to those of the defender. After a lengthy deliberation, they decided by a two-thirds majority:

\

$q$ Untel is guilty.

\

\noi When the time came for sentencing, at the general surprise, they decided, by a two-thirds majority:

\

$r$ Untel shall not be punished.

\

\noi How can this be? Just look at a miniature of Condorcet situations. One third of the Council voted yes for $p$ and $q$, but against $r$, quite logically. A third voted yes for $p$, against $q$, so yes for $r$ in all consistency. The last third voted against $p$ then yes for $q$ and $r$. {\bf Whence the paradoxical result!}

\

\noi My young listeners resented (that was obvious) the idea that the democratic rule of majority can lead (sometimes) to such inconsistencies, to aberrations where those who themselves voted the law by a majority could also violate it by a majority. Their disarray naturally led them to refuse the idea that this could be ineluctable and they sought actively solutions to this paradox. Of course, among those solutions there was one where it would suffice to revise the passed law, to repeal it. But, if none of the Elders changed his mind, no repeal would be possible since the vote would be reproduced just alike! One of them ended up with a more \lq\lq elaborate" solution. When we are faced with three issues where the answers yes to any two of them logically implies the answer yes to the third, it would be enough to examine only TWO of the issues out of three! This would avoid contradiction and {\it dissolve} the paradox. I made him observe that still remained a {\it subsidiary} question to settle: How to choose THE issue to discard among the three? Should we or can we consider the chronological order? for instance, or resort to another method? Which would remain to be defined.

\

\su{5 An apologue}  A City Council discusses the opportunity to build a community school and the place where it should be, possibly, built, knowing that there are only two available grounds, $U$ and $V$ . In short, let us say that the Council has to answer the following three questions:

\

$p$ \ Should we build a school?

$q$ \ Should we build a school on ground U?

$r$ \ Should we build a school on ground V?

\

\noi After a long deliberation, the Council answered {\bf yes} question $p$, by a majority vote. Then, with the same impetus, the Council answered {\bf no} question $q$, by a majority vote. The Mayor then said that a school will thus be constructed on ground $V$ . One of the counselors asks for a vote on question $r$, which excites the mirth of the assembly. He insists. Question $r$ is put to vote and, to their surprise, the Council's answer is {\bf no}, by a majority. We already know how this is possible, if the Council is divided into three equivalents groups, $H, K, L$, whose points of view are summed up in the following table:

\

$H \ \quad  p^+ \quad q^+\quad  r^-$

\

$K \ \quad  p^+ \quad q^- \quad r^+$

\

$L \ \quad  p^- \quad q^- \quad  r^-$

\

\noi Thus, after having raised issue $p$, depending on whether the next issue raised is $q$ or $r$, one can believe (in good faith) or argue (in bad faith) that the decision is in favor of ground $V$ or ground $U$, respectively. The order in which the questions are asked, if one decided to stick to the first two questions, would lead to three radically different decisions, according to cases. Thus, starting with $p$ and $q$, the Council would be led to decide to build a school on ground $V$ . If the first two addressed questions were $p$ and $r$, it would be led to construct a school on ground $U$. Finally, if $q$ and $r$ were put to a vote, the decision would be not to build a school at all. We hardly dare imagine the profit an unscrupulous Mayor, owner of ground $V$, could draw putting to vote questions $p$ and $q$, solely!

\

\noi The solution to have the Council vote in order to decide, itself, which of the three questions should be discarded is hardly satisfactory.

\

\noi The following three questions would have to be asked:

\

$P$ \ \ Should we discard question $p$?

$Q$ \ \ Should we discard question $q$?

$R$ \ \  Should we discard question $r$?

\

\noi We would face the following table which reflects the interests of the three groups, followed by the collective result of the Council's votes:

\

$H \ \  \  P^-  \ \ Q^+ \ \ R^-$

\

$K \ \  \ P^+ \ \  Q^- \ \ R^-$

\

$L \ \ \ P^- \ \  Q^- \ \ R^-$

\

\noi In other words, the Council decides that none of the three issues $p, q, r$ can be discarded. Unless, of course, they try once again to remove one of the three issues $P, Q, R$ ... We thus see, emerging, a regression from one Condorcet situation to another Condorcet situation, a kind of \lq\lq infinite descent" that would truly be infinite!

\

\noi It is true, all these examples have a somewhat \lq\lq schematic" aspect. However, if we do not always spot situations {\it \ˆ la Condorcet} with that sharpness, in every day life, it is probably because we often avoid to ask \lq\lq the third question!"

\

\su{6 A real situation}

\

\

Yet, here is a \lq\lq symptomatic" example, among others, I would like you to known.

\

\noi A poll was made in France, by Sofres, from 6 to 10 February 1987, with a national representative panel of 1,000 people, at the request of the weekly {\it Le Point}, about a possible reform shortening the presidential term (which was 7 years) to 5 years.

\

There were three options:

\

\noi $p$ \ It is desirable to make this reform before the presidential election of 1988.

\

\noi $q$ \  It is desirable to make this reform after the presidential election of 1988.

\

\noi $r$ \ It is better not to do this reform at all.

\

Here are the percentages: There are 10 percent with \lq\lq no opinion". Here is the distribution, among those who gave their opinions

\

39/90 \  chose \ $p$, 

25/90  \ chose \ $q$ 

26/90  \ chose \ $r$.

\

So 64/90 are for the reform, 51/90 do not want it to be made after the presidential election, and 65/90 do not want it done before!

\

Let the numbers talk for themselves!

[\lq\lq Laissons aux chiffres leur propre \Žloquence" !]

\

\noi This reform, however, ended up being done, much later (one easily understands why). Here is the time to add this: Of course, what we have said about the Condorcet paradox, hitherto, is about \lq\lq synchrony" (statics), the \lq\lq diachronic" study (dynamics), is much more complex, clearly (since it is necessary to take into account the evolution of individual choices and preferences, as well as any consultations leading to previous agreements, as in game theory).

\

An {\bf Infinite Society}, governed by an ultrafilter, would easily avoid all these disorders without having to resort to dictatorship ... as long as the ultrafilter is not trivial!

\

\head{\S2 Limits of families of sets}

\

\

\su{1 Preliminaries}

\

\

There is a special case when the upper and lower limits (see section {\bf 33} above) take a more set theoretic character than a topological one, when the space $E$ is equipped with its discrete topology. In this case, of course, the set $\Phi(E)$ of closed subsets of $E$ is none other than the set $\cal P(E)$ of all subsets of $E$, and the limits in this set are limits of a topology since a discrete space is locally compact. For a given family $\cal E = (E_i)_{i\in I}$  of subsets of $E$ and any filter $\cal F$ on $I$ whose grill is $\cal G$, by definition, we have
$$\underset{\cal F}\liminf \ \cal E = \underset{J \in \cal G} \bigcap \
\underset{i\in J}\bigcup \ E_i \ \text{and}\ \  \underset{\cal F}\limsup \ \cal E= \underset{J \in \cal F} \bigcap \ \underset{i\in J}\bigcup \ E_i $$
When the $E_i$'s are {\it just} sets, there is no need to specify the set $E$ of which they are subsets. The definitions remain unchanged whether $E$ is supposed to be the union of all these $E_i$ or any other set containing them.

\

\noi Let us insist on the following point. In the \lq\lq set-theoretic case", the upper and lower limits can also be written as follows:
$$\underset{\cal F}\liminf \ \cal E = \underset{J \in \cal F} \bigcup \
\underset{i\in J}\bigcap \ E_i \ \text{and}\ \  \underset{\cal F}\limsup \ \cal E= \underset{J \in \cal G} \bigcup \ \underset{i\in J}\bigcap \ E_i $$
This can easily be seen introducing, as in the \lq\lq topological" case, for each $x$, the set $I(x)=\{i : i\in I \ \text{and} \ x\in E_i\}$. Then  verify that we have:
$$x\in \underset{\cal F}\liminf \ \cal E \iff I(x)\in \cal F \iff x\in \underset{J \in \cal F} \bigcup \
\underset{i\in J}\bigcap \ E_i,$$
$$x\in \underset{\cal F}\limsup  \ \cal E \iff I(x)\in \cal G \iff x\in \underset{J \in \cal G} \bigcup \
\underset{i\in J}\bigcap \ E_i.$$

\

\su{2 Limits of sequences of sets}

\

\

The upper and lower limits of {\bf sequences} of sets appear in questions about measure and probability. They represent, in a sense, \lq\lq tails" of events and, thus, generalize the operations of union and intersection. When the filter $\cal F$ has a {\it countable} base, these two limits belong necessarily to the {\bf tribe} of measurable sets.

\

\noi Let $\cal E = (E_n)_{n\in \N}$  be a sequence of arbitrary sets and $\cal F$ be a filter on the set $\N$. When the sequence is monotone (increasing or decreasing) its limit relative to Fr\Žchet's filter [the filter of cofinite subsets in $\N$] always exists and is, depending on cases, either their union or their intersection. When the filter $\cal F$ has a countable base $\cal B$, the upper and lower limits both belong to the $\sigma$-algebra generated by the sets $E_n$. Indeed, for each base $\cal B$ of the filter $\cal F$, we have
$$\underset{\cal F}\liminf \ \cal E = \underset{J \in \cal B} \bigcup \
\underset{n\in J}\bigcap \ E_n \ \text{and}\ \  \underset{\cal F}\limsup \ \cal E= \underset{J \in \cal B} \bigcap \ \underset{n\in J}\bigcup \ E_n.$$
Now, I would like to give you an illustration of the use that can be made of set limits in a completely different domain.

\

\su{3 A preparatory lemma}

\

\

Let $\cal E = (E_i)_{i\in I}$ be a family of sets. For each ordered pair of sets $F$ and $M$, set $I[F, M] = \{i : i \in I \ \text{and} \  F \cap M = F \cap E_i\}$. Let then $\cal U$ be an ultrafilter on $I$, then set $L = \lim_\cal U \ \cal E$.

\

{\bf Lemma.}- {\it For each {\bf finite} set $F$, we have $I[F,L] \in \cal U$.}

\

\noi {\bf Indeed}, by definition, $i \in I[F, L]$ means that, for each $x \in  F \cap L$, we have $i \in I(x)$ and, for each $x \notin  F \cap L$, we have $i \notin  I(x)$. Therefore, we have
$$I[F,L]= \underset{x\in F\cap L} \bigcap \ I(x) \cap \underset{x\in F\stm L} \bigcap \ I \stm I(x).$$ 
Now, $x\in L \iff I(x)\in \cal U$ and  $x\notin L\iff I\stm I(x)\in \cal U$. So that the set $I[F,L]$ is the intersection of a finite number of subsets which belong to $\cal U$, whence the result. \qed

\

\su{4 \lq\lq Diagonals"}

\

\

Let $\cal E=(E_i)_{i\in I}$ be a family of sets. I shall say that a set $D$ is a {\bf diagonal} of the family $\cal E$ whenever, for each finite set $F$, the set of indices, $I[F, D]$, has the same cardinal as $I$, i.e., $|I[F, D]| = |I|$.

\

\noi In particular, a diagonal $D$ of the family $\cal E$ is thus a set of which every finite piece is a finite piece of at least one of the $E_i$'s, and still better, is the same finite piece from a number $|I|$ of the $E_i$'s. We can (pictorially) say that $D$ is thus {\bf upholstered} with finite pieces from the $E_i$'s.

\

\noi A filter $\cal F$  on a set $I$ is said to be {\bf  uniform} when each $X \in \cal F$ has the same cardinal as $I$, i.e., $|X| = |I|$. For that to be the case, it is necessary and sufficient that no subset $Y \inc I$ such that $|Y | < |I|$ belongs to $\cal F$. Set $\cal F_I =\{X : X\inc I \ \text{and} \  |I\stm X|<|I|\}$. This is a filter on I, the analog of Fr\Žchet's filter on $\N$ (except if $I$ is finite and not a singleton). \{When $I$ is finite, the set $\cal F_I$ is never a filter except if $I$ is a singleton.\} The uniform ultrafilters on $I$ are thus the ultrafilters which are finer than the filter $\cal F_I$. When $I$ is infinite, the uniform ultrafilters are not trivial. When $I$ is finite, there are no uniform ultrafilters on $I$ except in the trivial case when $I$ is a singleton. On a countable infinite set, the uniform ultrafilters are, precisely, the non trivial ultrafilters.

\

The preliminary lemma yields, immediately, the following result.

\

\su{5 Theorem} {\it For each family of sets, $\cal E = (E_i)_{i\in I}$, and each uniform ultrafilter $\cal U$ on $I$, the limit $\lim_\cal U  \cal E$ is a diagonal of the family $\cal E$.}

\

It can be shown, conversely, that each diagonal of the family $\cal E$ is a
limit of this family relative to a suitable uniform ultrafilter on $I$. 

\

How can diagonals be used?

\

\su{6 Applications}

\

\

I shall take a first example in Number theory.

\

{\bf The Semigroup N.}- Take an infinite subset $M \inc  \N$ and a family of subsets $\cal B = (B_m)_{m\in M}$  such that, for each $m \in M$, we have the  following: $B_m \inc  [0, m] \inc  B_m + B_m$. The \lq\lq additist arithmeticians" would call $B_m$ a base of the interval $[0, m]$. We know that diagonals of the family $\cal B$ exist. For each diagonal $D$ of $\cal B$, we have $D + D = \N$, that is $D$ is a base of $\N$.

\

\noi {\bf Indeed},  take any $n \in \N$. We will show that we have  $n \in D+D$ . From the hypothesis, we know that there exists an infinity of indices $m \in  M$  such that $[0,n]\cap D = [0,n]\cap B_m$. Taking such an $m \geqs n$ , we therefore have $n \in ([0,n]\cap B_m)+ ([0,n])\cap B_m)$ and, a fortiori, $n \in D+D$, as aforesaid. \qed

\

\noi However, we can do even better, in this case. For each subset $A \inc  \N$ and each integer $n$, it is customary to consider the number $r(A, n) = |\{(x,y) : x \in A, y \in  A, x + y = n\}|$ of representations of the integer $n$ as a sum of two integers from $A$ (taking into account the order of the terms in the sum). We then set $s(A) = \sup_{n\in \N} \ r(A, n)$.

\

\noi The result about the diagonal $D$ of the family $\cal B$ can be completed as follows. If we have $s(B_m)\leqs s$ for each $m\in M$ we also have $s(D)\leqs s$.

\

\noi {\bf Indeed}, taking again $n \in \N$ and the index $m \geqs  n$ such that $[0, n]\cap D = [0,n]\cap B_m$, observe that $r(D,n) = r([0,n]\cap D, n) \leqs r(B_m,n) \leq  s$ therefore $s(D) \leqs s$, as aforesaid. \qed

\

\noi This result appears in a paper by {\bf G. Grekos, L. Haddad, C. Helou, J. Pihko} [On the Erd\šs-Tur\'an conjecture, {\it  J. Number Theory }{\bf 102} (2003), $n^\circ$ 2, 339-352] with a different proof. It establishes the equivalence of two forms of an Erd\šs-Tur\'an conjecture, one {\it local} (strong), the other {\it global} (weak), so to say. See the 
article for more details.

\

The result easily generalizes as follows.

\

{\bf  Special semigroups.}- Let $S$ be an infinite commutative semigroup. Assume there is a special covering $(A_x)_{x\in S}$ in the following sense: For each $x\in S$, the set $A_x$ is a finite subset of $S$ such that $x\in A_x$ but $x\notin S+(S\stm A_x)$ and, moreover, $|S\stm \{y : y\in S, A_x \inc A_y\}|<|S|$ . Take a subset $M \inc  S$ such that $|M| = |S|$.

\

\noi For each $x\in M$, let $B_x$ be a base for $A_x$, i.e., $B_x \inc A_x \inc B_x + B_x$. Each diagonal $D$ of the family $B = (B_x)_{x\in M}$ is then a base for $S$, i.e., $D+D=S$. If, moreover, for each $x\in M$, we have $s(B_x)\leqs s$,we will also have $s(D) \leqs s$. The proof follows the same lines as in the preceding proof.

\

\noi {\bf Indeed}, take $x \in  S$ . From the hypothesis, the set $\{y \in M : A_x \cap  D = A_x \cap  B_y\}$ has the same cardinal as $M$ which is also the cardinal of $S$. Now, we have $|S\stm \{y : y \in S,  A_x \inc A_y\}| < |S|$. Therefore there exists at least one $y \in M$ such that, on the one hand, we have $A_x\cap D=A_x\cap B_y$ and, on the other, $A_x\inc A_y$. Since $B_y$ is a base for $A_y$, we have $x \in B_y +B_y$ hence $x \in (A_x \cap B_y)+(A_x \cap B_y)$ because $x \notin S + (S \stm  A_x)$. Since $A_x \cap D = A_x \cap B_y$, we have $x \in  D + D$ and, also, $r(D,x) \leqs  r(B_y,x) \leqs s(B_y) \leqs s$. \qed

\

This notion of diagonal appears naturally, at first, not through limits of sets, but in collective choices of a deliberative assembly or, if you will, as a \lq\lq nonstandard vision". Let us see, finally, how to express that and prove it.

\

\su{7 \lq\lq The genesis of diagonals"} Start with a family $B = (B_m)_{m\in M}$ with $M$ an infinite subset of $\N$, and each $B_m$ a base for the interval $[0,m]$. Suppose we have $s(B_m) \leqs s$ for each $m \in M$. Using the notations introduced by A. Robinson to denote \lq\lq hyperreal" objects, consider the family $^*\cal B = (B_m)_{m\in^*M}$. The set $M$ being infinite, there exists at least one hyperinteger $\mu \in ^* M \stm M$. \lq\lq By permanence", the subset $B_\mu$ is a base for the interval $[0,\mu]$ and we have $s(B_\mu) \leqs s$. So, for each hyperinteger $\nu \in [0,\mu]$, we have $\nu \in  B_\mu +B_\mu$. Let $D = B_\mu \cap \N$. For each integer $n \in \N$, we therefore have $n \in D + D$, which means that $D$ is a base for $\N$ and, moreover, we still have $s(D) \leqs s(B_\mu) \leqs s$.

\

So, to each hyperinteger $\mu\in ^*M \stm M$ corresponds a subset $D$ which is a diagonal of the family $\cal B$! I find that this way to look at things, which emerged first, is quite close to a certain intuition, by its simplicity. Each of the $B_m$'s being a base for $[0, m]$, what is more natural than go and see what happens when $m$ is an {\bf infinitely large} integer. That is how diagonals were born!

\

\head{A short bibliography}

\

1 ARROW K., {\sl Social choice and individual values}, John Wiley and Sons, New York, 1963.

\

2 BOURBAKI N., {\sl Fonctions d'une variable r\Želle}, FVR V.36, Appendice sur les corps de Hardy, Diffusion C.C.L.S., Paris, 1976; {\sl Topologie g\Žn\Žrale} TG I.43-46, Hermann, Paris, 1971.

\

3 CARTAN H., {\sl Th\Žorie des filtres}, C. R. Acad. Sc. Paris 205 (1937) 595-598; {\sl Filtres et ultrafiltres}, ibid., 777-779.
 
\

4 CHOQUET G., {\it Convergences}, Ann. Univ. Grenoble, Sect. Sci. Math. Phys. 23 (1947-1948) 57-112.

\

5 CONDORCET (Marquis de), {\sl Essai sur l'application de l'analyse \ˆ la probabilit\Ž des d\Žcisions rendues \ˆ la pluralit\Ž des voix}, Imprimerie Royale, Paris, 1785. (A photographic reprint has been published by Chelsea Publishing Company, 1972. The Discours has also been inserted in the following book: {\sl Condorcet, Sur les \Žlections et autres textes}, Corpus des \oe uvres de philosophie en langue fran\c caise, Fayard 1986.

\

6 GUILBAUD G.Th., {\sl Les th\Žories de l'int\Žr\t g\Žn\Žral et le probl\me logique de l'agr\Žgation}, \' Economie appliqu\Že, {\bf 5} (1952) $n^\circ$ 4, oct.-d\Žc., 501-551. (The paper was printed again as chapter II in the following book: {\sl \'El\Žments de la th\Žorie math\Žmatiques des jeux}, Monographies de recherches op\Žrationnelles, 9, Collection directed by G. Morlat, AFIRO, Dunod, Paris, 1968).

\

7 HADDAD L., {\sl Une repr\Žsentation des topologies}, C. R. Acad. Sc. Paris 255 (1962) 2702-2704; {\sl Sur quelques points de topologie g\Žn\Žrale. Th\Žorie des nasses et des tramails}, Ann. Fac. Sci. Univ. Clermont {\bf 44} fasc. 7 (1970) 3-80.

\

8 HADDAD L., {\sl Condorcet et les ultrafiltres}, in Math\Žmatiques finitaires et analyse nonstandard,  Publ. Math. Univ Paris VII, $n^\circ$ 31 tome 2 (1989) 343-360. Text of a lecture given at Luminy in 1985.

\

9 HADDAD L., {\'Elections, ultrafiltres, infinit\Žsimaux ou le paradoxe de Condorcet, in} Condorcet, Math\Žmaticien, \'Economiste, Philosophe, (p. 87-91), Homme politique, Colloque international, Paris, juin 1988, \Ždition Minerve.

\

10 HADDAD L., {\sl La double ultrapuissance}, S\Žminaire d'Analyse, Universit\Ž de Clermont II (1987-1988) expos\Ž 24.

\

11 \L O\'S J., {\sl O matrycach logicznych}, Prace Wroc\l awskiegoTowarzys\-twa Naukowego, Wroc\l aw, 1949.

\

12 \L O\'S J., {\sl Quelques remarques, th\Žor\mes et probl\mes sur les classes d\Žfinissables d'alg\bres, in} Mathematical interpretation of formal systems, Amsterdam, (1955) 98-113.

\

13 ROBINSON A., Non Standard Analysis, North-Holland, Amsterdam, 1966.

\

14 ELSHOLTZ Christian and LIST Christian, {\sl A Simple Proof of Sen's Possibility Theorem on Majority Decisions}, Elemente der Mathematik, {\bf 60} (2005) $n^\circ$ 2, 45-56.

\

\

This is a translation into English of a paper written in French, published in Tatra Mountains Mathematical Publications, {\sl L'ultrafiltre, un outil incomparable}, Tatra Mt. Math. Publ. {\bf 31} (2005), 131-176.

\noi It was also posted as {\tt arXiv:math/0702587v1} [math.HO] 20 Feb 2007.

\

\enddocument